\documentclass[5p]{elsarticle}

\usepackage{tikz}
\usepackage{tikz}
\usetikzlibrary{arrows.meta,positioning}
\graphicspath{{images/}}
\usepackage{amsmath}
\usepackage{latexsym, amssymb}

\usepackage{graphicx}
\usepackage{amsmath, amsbsy}
\usepackage{amsopn, amstext}
\usepackage{ifpdf}
\usepackage[bookmarks=false,hidelinks]{hyperref}
\usepackage{cancel, color}
\usepackage{epstopdf}
\usepackage{array,booktabs}
\usepackage{balance}

\def\ttimes{\,\rotatebox[]{-90}{$\ltimes$}\,}

\def\J{{\bf 1}}

\DeclareMathOperator{\Span}{Span}
\DeclareMathOperator{\Col}{Col}

\DeclareMathOperator{\lcm}{lcm}

\def\cal{\mathcal}

\def\pa{\partial}

\def\ra{\rightarrow}

\def\lra{\leftrightarrow}

\def\w{\wedge}

\def\a{\alpha}
\def\b{\beta}
\def\d{\delta}

\def\D{\Delta}

\def\0{{\bf 0}}

\def\A{\langle A \rangle}

\newcommand{\R}{{\mathbb R}}

\newcommand{\Z}{{\mathbb Z}}
\newcommand{\F}{{\mathbb F}}

\def\dsum{\mathop{\sum}\limits}

\newtheorem{theorem}{Theorem}[section]
\newtheorem{thm}[theorem]{Theorem}
\newtheorem{dfn}[theorem]{Definition}
\newtheorem{prp}[theorem]{Proposition}
\newtheorem{exa}[theorem]{Example}

\newtheorem{lem}[theorem]{Example}

\newtheorem{remark}{Remark}
\newtheorem{rem}[remark]{Remark}
\usepackage{algorithm}
\usepackage{algpseudocode}  

\newcommand{\placeholdergraphic}[1]{%
	\fbox{%
		\parbox[c][0.34\linewidth][c]{0.95\linewidth}{%
			\centering
			Placeholder figure\\[4pt]
			{\scriptsize image file unavailable}%
		}%
	}%
}

\begin{document}

\begin{frontmatter}
\title{\textbf{On Dimension-Varying Control Systems:\\A Universal State Space Approach}}


\author[tsinghua]{Feng Liu}\ead{lfeng@tsinghua.edu.cn}
\author[cas]{Daizhan Cheng}\ead{dcheng@iss.ac.cn}
\author[kth]{Xiaoming Hu}\ead{hu@kth.se}
\author[su]{Tielong Shen}\ead{tetu-sin@sopha.ac.jp}

\address[tsinghua]{State Key Laboratory of Power System Operation and Control, Dept. of Electrical Engineering, Tsinghua University, China}
\address[cas]{Key Laboratory of Systems and Control, Academy of Mathematics and Systems Sciences, Chinese Academy of Sciences, China}
\address[kth] {Department of Mathematics, KTH Royal Institute of Technology, Sweden}
\address[su]{Department of Engineering and Applied Sciences, Sophia University, Tokyo, Japan}

\begin{keyword} 
Dimension-varying (control) system, cross-dimensional Euclidean space, stabilization, disturbance decoupling, cross-dimensional projection.
\end{keyword}

\begin{abstract}               This paper develops a unified framework for the analysis and design of dimension-varying control systems by constructing an intrinsic quotient state space, $\Omega$. A significant challenge in non-fixed-dimensional systems is the lack of a common metric space that enables comparison of states across dimensions without relying on arbitrary external embeddings. To address this, we propose a cross-dimensional pseudo-metric $d_{\mathcal{V}}$ on $\mathbb{R}^{\infty}$ and derive $\Omega$ by identifying zero-distance representatives. We demonstrate that $\Omega$ preserves the essential topological and metric geometry of Euclidean space, providing the necessary foundation to extend fundamental control notions to the dimension-varying case. Specifically, we establish conditions for controllability, observability and stabilizability, and we address the complexities of Lipschitz switching and disturbance decoupling within this common space. The framework is further extended to hierarchical dimension-varying networks. The practical utility of the results is illustrated through a generator-removal-and-reconnection scenario in a three-machine power system. This case study demonstrates the use of translated representatives and projection/lift benchmarks, quantifies event-wise $d_{\mathcal{V}}$-gaps, and provides a finite-schedule dwell-time consistency check to validate the system’s structural transitions.
\end{abstract}

\end{frontmatter}

\section{Introduction}

Dimension-varying dynamic and control systems arise widely in natural and engineered systems. For example, on the Internet and in time-varying networks such as social, logistics, and power-supply networks, users join and leave frequently \cite{cas24}. In a biological system, cells are produced and die over time \cite{cha10}. Some man-made mechanical systems also vary in dimensions. For instance, the docking and undocking of spacecraft \cite{jia14,yan14}; the connecting and disconnecting of vehicle clutch systems while speed changes \cite{che20a}. The takeoff and landing of an aircraft can also be considered as a dimension-varying process. A similar situation arises in distributed power systems \cite{fen02} or electric vehicles with transmissions \cite{pak17}, just to mention a few.

To deal with such systems, various efficient techniques have been developed. For example, (1) using a global state space common to the different dynamics, together with reset maps\cite{hes02}; (2) gluing the state spaces of different modes together; \cite{jon14} provides a gluing condition for two state spaces of switched linear systems; (3) multi-model multi-dimensional ($M^3D$) system approach, which was initiated by \cite{ver06}.

Roughly speaking, the $M^3D$ approach is essentially to ``connect" the piecewise smooth trajectories together or measure the switching jump quantitatively.
\cite{ver13} proposed a conservative law to construct pseudo-continuous trajectories.
\cite{xue22} proposed a ${\cal K}{\cal L}$ function to estimate the instant jump disconnected trajectories.
Energy-based functions, such as a Hamiltonian function \cite{pak17} and a passivity-preserving energy \cite{gio25}, are used to measure switching discontinuities.

Another interesting phenomenon that motivates this paper is that a geometric object or a complex system may be described by models of different dimensions. For instance, in power systems, a single generator can be modeled as a $2$-, $3$-, $5$-, $6$-, or even $7$-dimensional dynamic system \cite{mac97}. In contemporary physics, string theory posits that the dynamics of strings model spacetime. But this model may have dimension 4 (special relativity), 5 (Kaluza-Klein theory), 10 (type 1 string), 11 (M-theory) or even 26 (Bosonic model) \cite{kak99}. One observes from this phenomenon that two models with different dimensions might be very similar or even equivalent.

Traditional opinion treats Euclidean spaces of different dimensions as completely separated or disconnected. If we consider
$$
\R^{\infty}:=\bigcup_{n=1}^{\infty}\R^n.
$$
Then each $\R^n$ is considered a component of $\R^{\infty}$. Topologically, each $\R^n$ is considered as a clopen (closed and open) subset of $\R^{\infty}$. (In this paper, such a topology is called a natural topology and denoted by ${\cal T}_n$.) The above phenomenon provides us a revelation that Euclidean spaces of different dimensions might be connected at certain points, such that the continuous functions, the vector fields and then the dynamic systems on different Euclidean spaces might be equivalent.

Based on this observation, this paper proposes a cross-dimensional distance $d_{{\cal V}}$ on $\R^{\infty}$ and identifies points with zero $d_{{\cal V}}$-distance, leading to the quotient space $\Omega=\R^{\infty}/\lra$. The point is not the use of a quotient construction alone, but the dimension-crossing quotient construction induced by $d_{{\cal V}}$: states belonging to Euclidean spaces of different dimensions become comparable elements of a single state space. This approach combines the aforementioned three techniques in a unified way. (1) The quotient space $\Omega$ serves as a common state space for dimension-varying systems, without imposing an external fixed-dimensional embedding as part of the state-space construction. Physical switching events, however, may still be represented by mode-dependent jump maps. (2) The equivalence relation $x\lra y$ if and only if $d_{{\cal V}}(x,y)=0$ glues dimension-compatible representatives into one quotient point and makes the induced space path connected. (3) The same distance $d_{{\cal V}}$ provides a quantitative measurement of dimension-changing switching gaps. Many switching descriptions in the literature can be formulated or tested under such a metric viewpoint, for example, a Hamiltonian function in \cite{pak17}, passivity energy in \cite{gio25}, and ${\cal K}{\cal L}$-class function in \cite{xue22}.

The central contribution of this paper is therefore to use this dimension-crossing quotient construction as a universal state space for dimension-varying systems. Since each fixed-dimensional component $\Omega^n$ is homeomorphic to the corresponding Euclidean space $\R^n$, and since the restriction of $d_{{\cal V}}$ to a fixed dimension is proportional to the ordinary Euclidean distance, $\Omega$ retains the finite-dimensional topology and metric geometry needed to reformulate classical control notions. Together with the lattice structure on $\Omega$, this allows dimension-varying trajectories, switching gaps, and control-design problems to be studied within the same quotient state space. The paper then discusses stability, controllability and observability, stabilization, disturbance decoupling, and large-scale hierarchical dimension-varying networks under this framework.

The rest of this paper is organized as follows: Section 2 provides mathematical preliminaries, including cross-dimensional Euclidean space and the dimension-keeping semi-tensor product of matrices. Section 3 formulates dimension-varying systems on the universal state space $\Omega$ and discusses cross-dimensional switching and stability. Section 4 introduces the differential-geometric structure of $\Omega$. Section 5 considers basic control-design problems for dimension-varying control systems. Section 6 discusses equivalent fixed-dimensional switching representations. Section 7 considers hierarchical networks with time-varying nodes. Section 8 presents an illustrative power-system example that instantiates the quotient-space construction. Section 9 gives brief concluding remarks.

Before ending this section, a list of notation is presented.

\begin{itemize}

\item $\R$: set of real numbers.

\item $\Z_+$: set of positive integers.

\item $\lcm(a,b)$ (or $a\vee b$): least common multiple of $a$ and $b$.

\item $gcd(a,b)$ (or $a\wedge b$): greatest common divisor of $a$ and $b$.

\item $\sup(a,b)$: least common upper bound of $a$ and $b$;

\item $\inf(a,b)$: greatest common lower bound of $a$ and $b$.

\item ${\cal M}_{m\times n}$: set of $m\times n$ real matrices.

\item$\otimes$: Kronecker product of matrices.

\item $\R^{\infty}:=\bigcup_{n=1}^{\infty}\R^{n}$.

\item $x\lra y$: $x,~y\in \R^{\infty}$ are equivalent.

\item $\Omega=\R^{\infty}/\lra $: the quotient space.

\item $\langle \cdot,\cdot\rangle_{{\cal V}}$: cross-dimensional inner product on $\R^{\infty}$,
(applicable to $\Omega$).

\item $\|\cdot\|_{{\cal V}}$: norm on $\R^{\infty}$.

\item $d_{{\cal V}}$: pseudo-metric on $\R^{\infty}$; it induces a metric on $\Omega$.

\item $m|n$: $m$ is a factor of $n$, i.e., $n/m$ is an integer.

\item $\J_n$: ${\underbrace{[1,\cdots,1]}_n}^T$.

\item $\ttimes$: dimension-keeping semi-tensor product.

\item $\d_n^i$: $\d_n^i=\Col_i(I_n)$.

\item $\d(t)$: Dirac function.

\item $\cong$: homeomorphism (of topological spaces).
\end{itemize}

\vskip 2mm

\section{Mathematical Preliminaries}

This section is based on \cite{che19,che19b,chepr}.

\subsection{Cross-Dimensional Euclidean Space}

Define a cross-dimensional Euclidean space as
$$
\R^{\infty}=\bigcup_{n=1}^{\infty}\R^n.
$$
It provides a meaningful state space for dimension-varying dynamic systems.

First, we need to turn it into a path-connected topological space, which allows the trajectories of a dimension-varying system to evolve continuously over such state space.

We define semi-tensor addition on $\R^{\infty}$.

\begin{dfn}\label{d.2.1.1} Let $x\in \R^m\subset \R^{\infty}$, $y\in \R^n\subset \R^{\infty}$, and $t=\lcm(m,n)$. Define
\begin{align}\label{2.1.1}
x\pm y:=\left(x\otimes \J_{t/m}\right)\pm \left(y\otimes \J_{t/n}\right).
\end{align}
\end{dfn}

Note that under this addition and the conventional scalar product of $x\in\R^{\infty}$ with $a\in \R$, $\R^{\infty}$ becomes a pseudo vector space, which is almost a vector space except that the zero element is not unique.

\begin{dfn}\label{d.2.1.2} Consider $\R^{\infty}$, assume $x\in \R^m\subset \R^{\infty}$, $y\in \R^n\subset \R^{\infty}$, and $t=\lcm(m,n)$.
\begin{itemize}
\item[(i)] The inner product of $x$ and $y$ is defined as
\begin{align}\label{2.1.2}
\langle x, y\rangle_{{\cal V}}:=\frac{1}{t}\langle x\otimes \J_{t/m}, y\otimes \J_{t/n}\rangle,
\end{align}
where $\langle \cdot,\cdot\rangle$ is the conventional inner product on $\R^t$.

\item[(ii)] The norm of $x$ is defined as
\begin{align}\label{2.1.3}
\|x\|_{{\cal V}}:=\sqrt{\langle x,x\rangle_{{\cal V}}}.
\end{align}

\item[(iii)] The cross-dimensional pseudo-distance between $x$ and $y$ is defined as
\begin{align}\label{2.1.4}
d_{{\cal V}}( x, y):=\|x-y\|_{{\cal V}}.
\end{align}
\end{itemize}

\end{dfn}

\begin{rem}\label{r.2.1.3}
\begin{itemize}
\item[(i)] Let $x,y\in \F^{\infty}$. $x$ and $y$ are said to be equivalent, denoted by $x\lra y$, if $d_{{\cal V}}(x,y)=0$.

\item[(ii)] It is easy to see that $x\lra y$, if and only if, there exist $\J_p$ and $\J_q$ such that
$$
x\otimes \J_p=y\otimes \J_q.
$$

\item[(iii)] Since distinct representatives may have zero $d_{{\cal V}}$-distance, $d_{{\cal V}}$ is a pseudo-metric on $\R^{\infty}$. Identifying zero-distance representatives gives the quotient space
$$
\Omega:=\R^{\infty}/\lra.
$$
The elements in $\Omega$ are denoted by
$$
\bar{x}=\{y\;|\;y\lra x\}.
$$

\item[(iv)]
Define
$$
\begin{array}{l}
\bar{x}\pm\bar{y}:=\overline{x\pm y},\\
\langle \bar{x},\bar{y}\rangle_{{\cal V}}:= \langle x, y\rangle_{{\cal V}},\\
\|\bar{x}\|_{{\cal V}}:=\|x\|_{{\cal V}},\\
d_{{\cal V}}(\bar{x},\bar{y}):=d_{{\cal V}}(x,y).
\end{array}
$$
They are all properly defined. After quotienting, $d_{{\cal V}}(\bar{x},\bar{y})$ is a genuine metric on $\Omega$. Hence $\Omega$ is a well-defined topological vector space \cite{kel63}. Meanwhile, it is also an inner product space and hence a metric space. But it is not a Hilbert space, because it is not complete.
\end{itemize}

(We refer to \cite{che19,chepr} for details of all above claims.)

\end{rem}

Next, we define a cross-dimensional projection.

\begin{dfn}\label{d.2.1.4}
 Let $\xi\in \R^n$. The projection of $\xi$ on $\R^m$, denoted by $\pi^n_m(\xi)$, is defined as
\begin{align}\label{2.1.5}
\pi^n_m(\xi):=argmin_{x\in \R^m}\|\xi - x\|_{{\cal V}}.
\end{align}
\end{dfn}

\begin{prp}\label{p.2.1.5} Let $t=\lcm(m,n)$. Then
\begin{align}\label{2.1.6}
\pi^n_m(\xi)=\Pi^n_m\xi,
\end{align}
where
\begin{align}\label{2.1.7}
\Pi^n_m=\frac{m}{t}(I_m\otimes \J^T_{t/m}) (I_n\otimes \J_{t/n}).
\end{align}
Moreover,
\begin{align}\label{2.1.8}
 \langle \xi - x_0, x_0\rangle_{{\cal V}}=0.
\end{align}
\end{prp}
Fig\ref{Fig.2.1} shows this.

\begin{figure}
\begin{center}
\setlength{\unitlength}{1cm}
\begin{picture}(5,4)\thicklines
\put(0,0){\line(1,0){5}}
\put(0,0){\vector(4,3){4}}
\put(4,3){\vector(0,-1){3}}

\put(1.9,1.9){$\xi$}
\put(2.5,0.1){$x_0$}
\put(4.1,1){$\xi-x_0$}

\put(3.7,0.3){\line(1,0){0.3}}
\put(3.7,0.3){\line(0,-1){0.3}}

\end{picture}\label{Fig.2.1}
\end{center}
\centerline{Fig.2.1~ Projection}
\end{figure}

It is also easy to prove that Pythagoras' theorem holds. That is,
\begin{align}\label{2.1.9}
\|\xi\|^2_{{\cal V}}=\|\xi-x_0\|^2_{{\cal V}}+\|x_0\|^2_{{\cal V}}.
\end{align}

\subsection{Dimension-Keeping Semi-Tensor Product}

Denote by
$$
{\cal M}:=\bigcup_{m=1}^{\infty}\bigcup_{m=1}^{\infty}{\cal M}_{m\times n}.
$$

\begin{dfn}\label{d.2.2.1}

Let $M\in {\cal M}_{m\times n}$, $N\in {\cal M}_{p\times q}$, $t=\lcm(n,p)$. The dimension-keeping semi-tensor product (DK-STP) of $M$ and $N$ is defined as
\begin{align}\label{2.2.3}
M\ttimes N:=\frac{n}{t}\left(M\otimes \J^{\mathrm{T}}_{t/n}\right)\left(N\otimes \J_{t/p}\right)\in {\cal M}_{m\times q}.
\end{align}
\end{dfn}

\begin{rem}\label{r.2.2.2}
In general, DK-STP is defined as
$$
M\ttimes N:=\left(M\otimes \J^{\mathrm{T}}_{t/n}\right)\left(N\otimes \J_{t/p}\right)\in {\cal M}_{m\times q}.
$$
The (\ref{2.2.3}) is called the weighted DK-STP. Weighted DK-STP has similar properties as DK-STP. Throughout this paper, we use only weighted DK-STP; hence (\ref{2.2.3}) is simply defined as DK-STP.
\end{rem}

Some fundamental properties are listed as follows:

\begin{prp}\label{p.2.2.3}
\begin{itemize}
\item[(i)] If $A\in {\cal M}_{m\times n}$ and $B\in {\cal M}_{n\times p}$, then
$$
A\ttimes B=AB.
$$
That is, the DK-STP is a generalization of the conventional matrix product.

\item[(ii)] Let $A\in {\cal M}_{m\times n}$, $B\in {\cal M}_{p\times q}$, and $t=\lcm(n,p)$. Then

\begin{align}\label{2.2.4}
A\ttimes B=A\Psi_{n\times p}B,
\end{align}
where $\Psi_{n\times p}\in {\cal M}_{n\times p}$, called a bridge matrix, which is defined as
\begin{align}\label{2.2.5}
\Psi_{n\times p}:=\frac{n}{t}\left(I_n\otimes \J^{\mathrm{T}}_{t/n}\right)\left(I_p\otimes \J_{t/p}\right).
\end{align}

\item[(iii)] (Distributivity) Assume $A$ and $B$ are of the same dimensions, say $A,B\in {\cal M}_{m\times n}$, then
$$
\begin{array}{l}
(A+B)\ttimes C=A\ttimes C+B\ttimes C,\\
C\ttimes (A+B)=C\ttimes A+C\ttimes B.
\end{array}
$$

\item[(iv)] (Associativity) Let $A,B,C\in {\cal M}$. Then
$$
(A\ttimes B)\ttimes C=A\ttimes(B\ttimes C).
$$

\item[(iv)] Consider ${\cal M}_{m\times n}$, with conventional matrix addition $+$
and DK-STP, $({\cal M}_{m\times n},+,\ttimes)$ is a ring\footnote{We refer to any standard Algebra textbook for related algebraic concepts, say, \cite{hun74}}.
\end{itemize}
\end{prp}

\begin{rem}\label{r.2.2.4} Comparing (\ref{2.2.5}) with (\ref{2.1.7}), one sees easily that
\begin{align}\label{2.2.6}
\Psi_{n\times m}=\Pi^m_n.
\end{align}
This identity is useful because we define them independently through two completely different ways.
But it provides a geometric meaning for DK-STP. That is, let $A\in {\cal M}_{m\times n}$ and $B\in {\cal M}_{p\times q}$, then
\begin{align}\label{2.2.7}
A\ttimes B=A\Pi^p_n B=A\pi^p_n(B).
\end{align}
It means first $\ttimes$ maps each column of $B$ into $\R^n$ and then $A$ maps vectors $x\in \R^n$ to $\R^m$.
\end{rem}

\section{Dimension-Varying Dynamic Systems}

\subsection{State Space: $\Omega$}

\begin{dfn}\label{d.3.1.1}
Let $\Sigma_i$, $i\in [1,s]$ be a finite set of dynamic systems with state spaces $\R^{n_i}$ respectively.
Assume there is a switching mapping $\sigma:\R_+\ra [1,s]$, which is piecewise constant and right continuous; then a dimension-varying system, consisting of modes
$\Sigma_i$, is defined as follows:
\begin{align}\label{3.1.1}
\Sigma(t)=\Sigma_{\sigma(t)},\quad t\in [0,\infty).
\end{align}
Denote the dimension-varying system briefly by $\Sigma=\{\Sigma_i\subset \R^{n_i}\;|\;i\in [1,s]\}$.
\end{dfn}

To begin with, we need to find a proper state space for dimension-varying systems. A fundamental requirement is that the trajectory of such a system can be continuous. A natural candidate is $\R^{\infty}$ with its topology induced by the distance $d_{{\cal V}}$, which is $\Omega$. Can we merge $\Sigma$ into $\Omega$? Yes! The following proposition shows that.

\begin{prp}\label{p.3.1.2}
$$
\Omega^n:=\Omega\bigcap \R^n\cong \R^n. \footnote{Where we consider $\Omega$ as its homeomorfic space $(\R^{\infty}, {\cal T}_d)$.}
$$
\end{prp}

\noindent{\it Proof.}
Let $x\in \R^n$. Then $\bar{x}\bigcap \R^n=x$, which is unique. Since any two norms on a finite-dimensional normed space are topologically equivalent\cite{tay80}, the claim is obvious. This fact can also be seen from the fact that the two distances are proportional.
\hfill $\Box$

\begin{rem}\label{r.3.1.3} There are two different topologies on $\R^{\infty}$: One is the natural topology, denoted by ${\cal T}_n$, which considers each $R^n$ as its component (precisely speaking, as its clopen (closed and open) subset). Another one is induced by the distance $d_{{\cal V}}$, denoted by ${\cal T}_d$, which makes $\R^{\infty}$ (precisely speaking $\Omega$) a path-connected topological space. Figure \ref{Fig.3.1} shows these two topologies.\footnote{We refer to any Topology textbook for concepts in topology, say, \cite{kel55}.}
\end{rem}

\begin{figure}
\begin{center}
\setlength{\unitlength}{5.5mm}
\begin{picture}(9,15)(2,-5)
\thicklines
\put(8,3){\line(1,0){7}}
\put(8,3){\line(-1,-1){2}}
\put(13,1){\line(1,1){2}}
\put(13,1){\line(-1,0){7}}
{\color{red}
\put(8,2){\line(1,1){2}}
\put(8,2){\line(1,-1){2}}
\put(12,2){\line(-1,1){2}}
\put(12,2){\line(-1,-1){2}}
}
{
\put(7.8,0){\line(1,1){4}}
}
\put(5.5,0){\line(1,0){9}}
\put(5.5,0){\line(0,1){5}}
\put(5.5,5){\line(1,0){9}}
\put(5.5,5){\line(1,1){2}}
\put(6,4.2){$\Omega^{m\vee n}$}
\put(0,0){\line(1,0){3}}
\put(0,2){\line(1,1){1}}
\put(0,2){\line(1,0){3}}
\put(4,3){\line(-1,-1){1}}
\put(4,3){\line(-1,0){3}}
\put(1,6){\vector(-1,-1){1}}
\put(1,6){\vector(1,0){3}}
\put(1,6){\vector(0,1){1}}
\put(4,3){\line(-1,0){3}}
\put(4.2,0){$\R^1$}
\put(4.2,2.8){$\R^2$}
\put(4.2,6){$\R^s$}
\put(2,3.8){$\vdots$}
\put(13,2.2){$\Omega^n$}
\put(9.3,3.05){$\Omega^m$}
\put(11.5,4.2){$\Omega^{m\wedge n}$}
\put(1,8){$(\R^{\infty}, {\cal T}_n)$}
\put(8,8){$(\R^{\infty}, {\cal T}_d)\cong \Omega$}
\put(5,8){vs}
\end{picture}\label{Fig.3.1}
\end{center}
\vskip -2.5cm
\centerline{Fig.3.1~Comparing $\R^{\infty}$ with $\Omega$ }
\end{figure}

Proposition \ref{p.3.1.2} is of fundamental importance, because it leads to the following result.

\begin{thm}\label{t.3.1.4}
\begin{itemize}
\item[(i)] $\R^n$ with its classical Euclidean topology is homeomorphic to $\Omega^n=\R^n\cap \Omega$ with its subspace topology inherited from $\Omega$.
\item[(ii)]
Consider a sequence $\{\bar{x}_k\;|\; k=1,2,\cdots\}\subset \Omega$. Assume
$$
\lim_{k\ra \infty}\bar{x}_k=\bar{x}_0,
$$
and $x_0\in \bar{x}_0\bigcap \R^n$. Let
$$
x_{k_i}=\bar{x}_{k_i}\bigcup \R^n
$$
be a sequence in $\R^n$. Then
\begin{align}\label{3.1.101}
\lim_{i\ra \infty}x_{k_i}=x_0.
\end{align}
\end{itemize}
\end{thm}

\noindent{\it Proof.}
\begin{itemize}

\item[(i)] Let $x,y\in \R^n$. By definition, we have
$$
d_{{\cal V}}(\bar{x},\bar{y})=\frac{1}{\sqrt{n}}d(x,y).
$$
Note that over $\R^n$ the $\frac{1}{\sqrt{n}}$ is a constant number; then the topology induced by $d(x,y)$ is the same as the topology induced by $d_{{\cal V}}$. The former is the classical Euclidean topology of $\R^n$, and the latter is the subspace topology of $\R^n$ inherited from $\Omega$. The conclusion follows.

\item[(ii)] Since
\begin{align}\label{3.1.102}
d(x_{k_i},x_0)=\sqrt{n}d_{{\cal V}}(\bar{x}_{k_i},\bar{x})\xrightarrow{i\ra \infty} 0,
\end{align}
the conclusion is obvious.
\end{itemize}
\hfill $\Box$

\begin{rem}\label{r.3.1.401}
\begin{itemize}
\item[(i)] Fact (i) of Theorem \ref{t.3.1.4} is surprising. In fact, $\Omega$ is a quotient space of $(\R^{\infty}, {\cal T}_n)$. Hence the quotient topology on $\Omega$ is coarser than ${\cal T}_n$, and it is easy to verify that it is strictly coarser than ${\cal T}_n$. One may therefore expect the subspace topology on $\Omega^n\subset \Omega$ to be coarser than the classical topology on $\R^n$. The theorem shows that this does not happen: $\Omega^n$, with its subspace topology, is homeomorphic to ${\cal T}_n\vert_{\R^n}$, namely, the classical Euclidean topology. Moreover, the proof shows a stronger metric compatibility on each fixed-dimensional component, since for $x,y\in\R^n$,
$$
d_{{\cal V}}(\bar{x},\bar{y})=\frac{1}{\sqrt{n}}d(x,y).
$$
Thus the homeomorphism accounts for the inherited topology, while the restriction of $d_{{\cal V}}$ accounts for the inherited metric geometry. This distinction is important: the quotient construction preserves the finite-dimensional Euclidean structure locally on each $\Omega^n$, while the full space $\Omega$ provides the additional cross-dimensional comparability. This fact makes it convincing to take $\Omega$ as the state space of dimension-varying systems.

\item[(ii)] Fact (ii) shows that cross-dimensional convergence of a sequence $\{\bar{x}_n\}$ ensures the convergence of $x_{n_i}\in \bar{x}_{n_i}\bigcap \R^{n_i}$ for any $n_i\in \Z_+$. Roughly speaking, the convergence in $\Omega$ implies the convergence in each $\Omega^n\cong\R^n$.
\end{itemize}
\end{rem}

Summarizing the above arguments, one sees that a dimension-varying system $\Sigma$ can be merged into $\Omega$ by merging each mode $\Sigma_i$ into the corresponding component $\Omega^{n_i}$. Then the dynamics of each subsystem evolving on $\Omega^{n_i}$ is represented without changing its local fixed-dimensional dynamics on $\R^{n_i}$.

The preceding result should be read in two layers. The homeomorphism $\Omega^n\cong\R^n$ preserves the fixed-dimensional topology, while the proportional relation between $d_{{\cal V}}$ and the Euclidean distance on $\R^n$ preserves the metric geometry on that component. Across different dimensions, the inner product and distance of $\Omega$ compare representatives after lifting them to a common multiple dimension. The following example depicts this quotient-metric geometry.

\begin{exa}\label{e.3.401} Consider a triangle $\bigtriangleup ABC\in \Omega$ (see Figure \ref{Fig.3.2}), where $A,~B,~C\in \Omega$, say, $A=(1,2,-1)^{\mathrm{T}}\in \Omega^3$, $B=(2,0,-1,3)^{\mathrm{T}}\in \Omega^4$, and $C=(1,2,-1,-2,1)^{\mathrm{T}}\in \Omega^5$.

\begin{figure}
\begin{center}
\setlength{\unitlength}{12mm}
\begin{picture}(7,4)
\thicklines
\put(6,1){\vector(-1,0){5}}
\put(1,1){\vector(2,1){4}}
\put(5,3){\vector(1,-2){1}}
\put(5,3.2){$A$}
\put(0.5,0.8){$B$}
\put(6.2,0.8){$C$}
\put(3.5,0.5){$a$}
\put(2.8,2.2){$c$}
\put(5.7,2){$b$}
\end{picture}\label{Fig.3.2}
\end{center}

\centerline{Fig.3.2~A Triangle in $\Omega$ }
\end{figure}

Set three sides as
$$
cb:=B-C,~ac:=C-A,~ba:=A-B.
$$
Then it is easy to calculate that
$$
\begin{array}{l}
a=\|cb\|_{{\cal V}}=1.7607,\\
b=\|ac\|_{{\cal V}}=1.9833,\\
c=\|ba\|_{{\cal V}}=2.5499.
\end{array}
$$
$$
\begin{array}{l}
\cos{\angle A}=\frac{-\langle ba,cb\rangle_{{\cal V}}}{b c}=0.7252,\\
\cos{\angle B}=\frac{-\langle cb,ba\rangle_{{\cal V}}}{a c}=0.6312,\\
\cos{\angle C}=\frac{-\langle cb,ac\rangle_{{\cal V}}}{a b}=0.0764.\\
\end{array}
$$
It follows that
$$
\begin{array}{l}
\angle A=43.5177^o,\\
\angle B=50.8619^o,\\
\angle C=85.6200^o.
\end{array}
$$
It is straightforward to verify the Law of Sines:
$$
\frac{a}{\sin{A}}=\frac{b}{\sin{B}}=\frac{c}{\sin{C}}=2.5571.
$$
It is also easy to verify that the Law of Cosines, Law of Tangents, etc. are also correct.

In fact, these calculations use the inner product and norm induced by $d_{{\cal V}}$ after the involved representatives are compared in a common dimension. Therefore, the usual metric identities of Euclidean inner-product geometry, such as the laws of sines and cosines, remain valid for such quotient-metric calculations in $\Omega$.

This observation should not be interpreted as saying that the whole $\Omega$ is a single finite-dimensional Euclidean space. Rather, each fixed-dimensional component inherits Euclidean topology and metric geometry, while the quotient and lattice structures provide cross-dimensional comparability.
\end{exa}

\subsection{Lattice Structure on $\Omega$}

\begin{dfn}\label{d.3.15.1} \cite{abb69,bur81} Let ${\cal L}\neq \emptyset$ be a partially ordered set with the order $\prec$. If for any two elements $a,~b\in {\cal L}$ there is a least common upper bound, denoted by
$\sup(a,b)\in {\cal L}$, and a greatest common lower bound, denoted by $\inf(a,b)\in {\cal L}$, then ${\cal L}$ is called a lattice.
\end{dfn}

\begin{exa}\label{e.3.15.2} Consider the set of positive integers $\Z_+$.
Define an order
\begin{align}\label{3.15.1z}
a\prec b,\quad \mbox{iff}~ a|b,\quad a,b\in \Z_+.
\end{align}
It is easy to verify that
$$
\sup(a,b)=\lcm(a,b),\quad \inf(a,b)=\gcd(a,b).
$$
Hence $(\Z_+,\prec)$ is a lattice, called the multiplication-division (MD-1) lattice.
\end{exa}

\begin{dfn}\label{d.3.15.3} Consider
$$
\Omega=\bigcup_{n=1}^{\infty}\Omega^n.
$$
As aforementioned that $\Omega^n$ is homeomorphic to $\R^n$, $n=1,2,\cdots$.
$\Omega^m$ is said to be an ${\cal L}$-subspace of $\Omega^n$ and denoted by
\begin{align}\label{3.15.1}
\Omega^m\prec \Omega^n
\end{align}
if $m|n$.

\end{dfn}

\begin{prp}\label{p.3.15.4} Under the order defined by (\ref{3.15.1})
$$
{\cal L}_{\Omega}:=\left(\{\Omega^n\;|\;n=1,2,\cdots\},\prec \right)
$$
is a lattice, where
\begin{align}\label{3.15.3}
\begin{array}{l}
\sup(\Omega^m,\Omega^n)=\Omega^{m\vee n},\\
\inf(\Omega^m,\Omega^n)=\Omega^{m\wedge n}.
\end{array}
\end{align}
\end{prp}

\noindent{\it Proof.} Assume $\Omega^t\prec \Omega^m$ and $\Omega^t\prec \Omega^n$. By definition, $t|m$ and $t|n$. That is, $t$ is a common divisor of $m$ and $n$. Hence, $t|(m\wedge n)$. It follows that $\Omega^{m\wedge n}$ is the greatest common lower bound of $\Omega^m$ and $\Omega^n$. Similarly, we can prove $\Omega^{m\vee n}$ is the lowest common upper bound of $\Omega^m$ and $\Omega^n$.
\hfill $\Box$

It is obvious that ${\cal L}_{\Omega}:=\left(\{\Omega^n\;|\;n=1,2,\cdots\},\prec \right)$ is isomorphic to MD-1.

\begin{dfn}\label{d.3.15.5} Let ${\cal L}$ be a lattice.
\begin{itemize}
\item[(i)] Let ${\cal H}\subset {\cal L}$ be a subset. ${\cal H}$ is called a sub-lattice, if
$$
\begin{array}{l}
\sup(a,b)\in {\cal H},\\
\inf(a,b)\in {\cal H},\quad a,~b\in {\cal H}.
\end{array}
$$
\item[(ii)] A sub-lattice ${\cal H}$ is called an ideal if, whenever $a\in {\cal L}$ and $a\prec h$ for some $h\in {\cal H}$, then $a\in {\cal H}$.
\item[(iii)] A sub-lattice ${\cal H}$ is called a filter if, whenever $a\in {\cal L}$ and $h\prec a$ for some $h\in {\cal H}$, then $a\in {\cal H}$.

\item[(iv)] A filter ${\cal F}$ is called a principal filter if there exists a unique smallest element $x_0$ such that
$$
h\succ x_0, \quad h\in {\cal F}.
$$
The ${\cal F}$ is denoted by ${\cal F}(x_0)$.
\item[(v)] An ideal ${\cal I}$ is called a principal ideal if there exists a unique largest element $x_0$ such that
$$
h\prec x_0, \quad h\in {\cal H}.
$$
The ${\cal I}$ is denoted by ${\cal I}(x_0)$.
\end{itemize}
\end{dfn}

Consider $\bar{x}\in \Omega$. It is easy to see that there exists a unique $x_1\in \bar{x}$, such that
$$
\bar{x}=\{x_n=x_1\otimes\J_n\;|\; n=1,2,\cdots\}.
$$
Then we define
$$
\dim(\bar{x}):=\dim(x_1),\quad \bar{x}\in \Omega.
$$
\begin{exa}\label{e.3.15.6} Consider $\R^{\infty}$.
Assume $x\neq y\in \R^{\infty}$. A partial order $\prec$ is defined as follows: $x\prec y$, iff, there exists a $\J_s$ such that $x\otimes \J_s=y$. Then $\left(\R^{\infty},\prec\right)$ is a lattice.
\begin{itemize}
\item[(i)] $x$ and $y$ are comparable, iff, $x\lra y$.
\item[(ii)] Let $\bar{x}\in \Omega$ and $x_1\in \bar{x}$ be the smallest element. Then
$\bar{x}$ is a principal filter generated by $x_1$.
\end{itemize}
\end{exa}

\begin{dfn}\label{d.3.15.7} Assume $\bar{x}\in \Omega$ and $\dim(\bar{x})=n_0$. Then the tangent bundle of $\bar{x}$ is
\begin{align}\label{3.15.4}
T_{\bar{x}}=\bigcup_{k=1}^{\infty}T(\R^{kn_0}).
\end{align}
\end{dfn}

\subsection{Lipschitz Switching}

This subsection considers the switchings in a dimension-varying system.

\begin{dfn}\label{d.3.2.1} Consider a cross dimensional system $\Sigma=\{\Sigma_i\subset \R^{n_i}\;|\;i\in[1,s]\}$.
\begin{itemize}
\item[(i)]
Let
$t$ be a switching moment such that $x(t^-)\in \Sigma_p$ and $x(t^+)\in \Sigma_q$. The switching is called a Lipschitz switching with respect to $L^r_s>0$, if
\begin{align}\label{3.2.1}
\|x(t^+)\|_{{\cal V}}\leq L^p_q\|x(t^-)\|_{{\cal V}}.
\end{align}
\item[(ii)] The system $\Sigma$ is said to have Lipschitz switching, if there exists a common $L>0$ such that all switchings are Lipschitz with respect to $L$.
\end{itemize}
\end{dfn}

Through this paper we assume there exists a Lipschitz function $\varphi$, such that
\begin{align}\label{3.2.2}
\|x(t^+)\|_{{\cal V}}=\varphi(\|x(t^-)\|_{{\cal V}}).
\end{align}

In the following example, we propose some Lipschitz switches.

\begin{exa}\label{e.3.2.3}
\begin{itemize}
\item[(i)] Constant Linear Mapping:

Assume there exist a set of constant matrices $W^p_q\in {\cal M}_{n_q\times n_p}$ such that the switching $\sigma(t): x(t^-)\in \Sigma_p \ra x(t^+)\in \Sigma_q$ satisfies
$$
x(t^+)=W^p_qx(t^-).
$$

Let $A\in {\cal M}$, then the operator norm of $A$ on $\R^{\infty}$ is defined by
$$
\|A\|_{{\cal V}}:=\sup_{0\neq x\in \R^{\infty}}\frac{\|A\ttimes x\|_{{\cal V}}}{\|x\|_{{\cal x}}}.
$$
It has been proved in \cite{che19b} that (see also \cite{che23} Proposition 13.4.2)
if $A\in {\cal M}_{m\times n}$, then
\begin{align}\label{3.2.3}
\|A\|_{{\cal V}}=\sqrt{\frac{n}{m}\sigma_{\max}(A^TA)}.
\end{align}

Using (\ref{3.2.3}) we have that
\begin{align}\label{3.2.4}
\begin{array}{l}
 L^p_q=\|W^p_q\|_{{\cal V}}
=\sup_{0\neq x\in \R^{\infty}}\frac{\|W^p_qx\|_{{\cal V}}}{\|x\|_{{\cal V}}}\\
~\\
=\sqrt{\frac{n_p}{n_q}\sigma_{\max}((W^p_q)^{\mathrm{T}}W^p_q)}.
\end{array}
\end{align}

\item[(ii)] Nearest Jump:

Assume the trajectory jumps at $t$ from $p$-th mode to $q$-th mode. That is, from $\R^{n_p}$ to $\R^{n_q}$. It jumps to the nearest point. That is,
\begin{align}\label{3.2.5}
x(t^+)=\Pi^{n_p}_{n_q}x(t^-).
\end{align}

This is the particular case of (i) with $W^p_q=\Pi^{n_p}_{n_q}$.

\item[(iii)] Differentiable Jump:

In general, assume the jump functions
$$
x^q_i=x^q_i(x^p_1,\cdots,x^p_{n_p}),\quad i\in [1,n_q],
$$
are differentiable functions. Then
$$
W^p_q(t)=\frac{\pa(x^q_1,\cdots,x^q_{n_q})}{\pa(x^p_1,\cdots,x^p_{n_p})}(t).
$$
The Lipschitz condition ensures the existence of a uniform bound.

\end{itemize}
\end{exa}

\begin{rem}\label{r.3.2.4}
\begin{itemize}
\item[(i)] Nearest jump is a canonical geometric switching law under $d_{{\cal V}}$. In applications, it may be used either as an idealized switching rule or as a benchmark against which realized switching policies are compared.
\item[(ii)] If the transition occurs in a small neighborhood of the closest point, the mismatch can be treated as an error term. The subsequent discussion remains valid.
\item[(iii)] A common situation for the switches of a network is: some new recruits join the network, or some members leave the network.
If a switch is caused by some members leaving the network, it is obvious that the Lipschitz condition is satisfied. As new recruits join the network, the Lipschitz condition means they must meet certain requirements. Otherwise, they may destroy the network structures.
\item[(iv)] In practice, the jumping functions are mostly experimentally derived. For instance, in the vehicle clutch transient process \cite{tem18} provides such a jumping function.
\end{itemize}
\end{rem}

Define the switching gap by
\begin{align}\label{3.2.6}
\d t:=d_{{\cal V}}(x(t^-), x(t^+)).
\end{align}

By definition, we have the following result.

\begin{prp}\label{p.3.2.5} Consider $\Sigma=\{\Sigma_i\subset \R^{n_i}\;|\;i\in [1,s]\}$. Its trajectory is continuous at $t$, if and only if,
$x(t^-)\lra x(t^+)$, i.e., $d(t)=d_{{\cal V}}(x(t^-),x(t^+))=0$.
\end{prp}

Observing Figure \ref{Fig.3.1}, it is clear that a switching is continuous if and only if $x(t^-)\in \R^{n_p}\bigcap \R^{n_q}=\R^{n_k}$, where $n_k=\gcd(n_p,n_q)$.

\begin{exa}\label{e.3.2.7}
\begin{itemize}
\item[(i)] Assume $\Sigma_1\subset \R^n$ and $\Sigma_2\subset \R^m$, where $m<n$ and $\Sigma_2$ is obtained from $\Sigma_1$ by dropping $n-m$ agents. Without loss of generality, we assume the last $n-m$ agents are dropped.
 Then the transition matrix becomes
 \begin{align}\label{3.2.7}
 W^n_m=[I_m,0_{m\times (n-m)}].
 \end{align}
It is obviously a Lipschitz switching.
\item[(ii)] Assume $\Sigma_1\subset \R^n$ and $\Sigma_2\subset \R^m$, where $m>n$ and $\Sigma_2$ is obtained from $\Sigma_1$ by adding $n-m$ agents. We assume the addition is obtained by nearest jumping. Then the transition matrix
is
\begin{align}\label{3.2.8}
W^n_m=\begin{bmatrix}
I_n\\
\Pi^n_{m-n}\\
\end{bmatrix},
\end{align}
where $\Pi^n_{m-n}$ can be replaced by any matrix with full row rank, or even a (full row rank) Lipschitz mapping.
\item[(iii)] If a switch consists of dropping and adding, then
\begin{align}\label{3.2.9}
W^n_m=W^s_mW^n_s,\quad s<n,
\end{align}
where $W^n_s$ is of the form (\ref{3.2.7}) and $W^s_m$ is of the form (\ref{3.2.8}).
\end{itemize}
\end{exa}

\subsection{Stability of Dimension-Varying Systems}

\begin{dfn}\label{d.3.3.1} Consider a dimension-varying system $\Sigma=\{\Sigma_i\subset \R^{n_i}\;|\;i\in [1,s]\}$.
\begin{itemize}
\item[(i)] An $\bar{x}_0\in \Omega$ is called an equilibrium of $\Sigma$, if $x^i_0\in \Omega^{n_i}\bigcap \bar{x}_0$ $\forall i\in [1,s]$ exist and $x^i_0$ is an equilibrium of $\Sigma_i$, $i\in [1,s]$.
\item[(ii)] $\Sigma$ is said to be (globally) stable at $\bar{x}_0$, if $\Sigma_i$ is stable at $x^i_0$, $\forall i\in [1,s]$.
\item[(iii)] $\Sigma$ is said to be locally stable at $\bar{x}_0$, if $\Sigma_i$ is locally stable at $x^i_0$, $\forall i\in [1,s]$.
\end{itemize}
\end{dfn}

\begin{dfn}\label{d.3.3.2} Consider a dimension-varying system $\Sigma=\{\Sigma_i\subset \R^{n_i}\;|\;i\in [1,s]\}$.
Let the jumping times be $0<t_1<t_2<\cdots$.
\begin{itemize}
\item[(i)]
The minimum dwell time is defined as
$$
\D_m:=\inf_{i\in \Z_+}(t_{i+1}-t_i)\geq \epsilon>0.
$$
\item[(ii)] Let $P_i$ be the set of non switching time intervals $[t_i,t_{i+1})$, when the model $\Sigma_i$ is active. Hence for each non-switching period $[t_{j},t_{j+1})$ there is a unique $P_i$ such that
$$
[t_j,t_{j+1})\in P_i.
$$
\end{itemize}
\end{dfn}

\begin{dfn}\label{d.3.3.3} Consider a dimension-varying system $\Sigma=\{\Sigma_i\subset \R^{n_i}\;|\;i\in [1,s]\}$.
 Let $n=\lcm(n_1,\cdots,n_d)$ and ${\cal L}$ be a positive definite function on $\R^n$ with
$$
L(x)
\begin{cases}
=0,\quad x=x_0,\quad x_0\in \bar{x}_0\\
>0,\quad \mbox{Otherwise},
\end{cases}
$$
 $L(x)$ is called a Lyapunov function for $\Sigma$ if
\begin{align}\label{3.3.1}
\frac{d}{dt}L(x^i(t)\otimes \J_{n/n_i})<0,\quad x^i(t)\neq x^i_0,\quad t\in p_i.
 \end{align}
 \end{dfn}

Then we have the following result.

\begin{thm}\label{t.3.3.3} Consider a dimension-varying system $\Sigma=\{\Sigma_i\subset \R^{n_i}\;|\;i\in [1,s]\}$.
The system $\Sigma$ is stable if there exists a Lyapunov function, the switching gaps are uniformly bounded, and the minimum dwell time $\D_m$ is large enough.
\end{thm}

\noindent{\it Proof.} Using Theorem \ref{t.3.1.4}, one sees easily that the restriction $L|_{\Sigma_i}$ is a Lyapunov function. Hence $L$ is essentially a common Lyapunov function for all models. Now the problem is that the trajectories are not guaranteed to be continuous. Hence, we have a jump gap $\d_{t_k}$ at jump moment $t_k$. Assume the dwell time $\D_d$ is large enough. Then we have
$$
\|\bar{x}_{t^+_{k}}\|_{{\cal V}}-\|\bar{x}_{t^-_{k+1}}\|_{{\cal V}}> \d_{t_k}.
$$
As long as the $\d_{t_k}$, $k=1,2,\cdots$ are uniformly bounded, the Lyapunov function can be used to find a finite dwell time $\D_d>0$ such that the decrease in norm during the dwell time is larger than the increase in norm by the preceding jump.

It remains to prove the existence of a uniform bound for $\d_{t_k}$, $k=1,2,\cdots$.

Since the switchings are assumed to be Lipschitz, the required dwell time can be obtained.
\hfill $\Box$

\noindent\textit{Remark.} For an affine switching law $x(t^+)=W^p_qx(t^-)+\eta^p_q$, the offset $\eta^p_q$ contributes to $\d_{t_k}$ and should be included when verifying the uniform gap-bound assumption; it does not modify the Lyapunov decrease argument.

We give a simple example to describe this.

\begin{exa}\label{e.3.3.4} Consider a dimension-varying system $\Sigma=\{\Sigma_1\subset \R^2,\Sigma_2\subset \R^4\}$, where
\begin{align}\label{3.3.2}
\begin{array}{ll}
\Sigma_1:&\dot{x}^1(t)=\begin{bmatrix}
3&2\\
-10&-6
\end{bmatrix}x^1(t):=A_1x^1(t),\\
~&x^1\in \R^2.\\
\Sigma_2:&\dot{x}^2(t)=\begin{bmatrix}
-5&1&0&1\\
1&-3&0&1\\
-1&0&-2&0\\
0&1&0&-2\\
\end{bmatrix}x^2(t):=A_2x^2(t),\\
~&x^2\in \R^4.
\end{array}
\end{align}
Using (\ref{2.1.7}), we have
$$
\Pi^2_4x-x=x\otimes \J_2-x=0.
$$
So it causes no transition gap.
Using (\ref{3.2.3}) yields
$$
\|\Pi^4_2\|_{{\cal V}}=2.
$$
Then
$$
\|\Pi^4_2x-x\|_{{\cal V}}\leq (\|Pi^4_2\|_{{\cal V}}+1)\|x\|_{{\cal V}}=3\|x\|_{{\cal V}}.
$$
It is easy to verify that
$$
\|e^{3.6A_1}\|=0.3201.
$$
It is clear that when the dwell time $\D_d=3.6$, $\Sigma$ is stable.

Note that the norm estimation is equivalent to setting $L(x)=\|x\|$, where $\|\cdot\|$ is the standard norm on $R^4$.
\end{exa}

\subsection{Robustness}

It was emphasized in \cite{xue22} that in dimension-varying systems there are state impulses appearing at the switching moments. Under our model, it can be treated as a disturbance. We use linear switching modes to describe this.
Consider a dimension-varying system
\begin{align}\label{3.4.1}
\Sigma=\{\Sigma_k\subset \R^{n_k}\;|\;k\in [1,s]\},
\end{align}
where $\Sigma_k$ is described by
 \begin{align}\label{3.4.2}
\dot{x}(t)=A_k\ttimes (x(t)+\eta_k(t)w),\quad t\in [t_i,t_{i+1}), k\in [1,s],
\end{align}
where $\eta_k(t)\in \R^{\ell}$ is a disturbance of dimension $\ell$, which is uncertain but measurable.

If it is the state impulse, it can be described as
$$
\eta_k(t)=\d(t)gv,
$$
where $\d(t)$ is a Dirac function, $g\in \R^{\ell}$ is of unit norm, where typically $\ell$ is in between $\dim(x(t^-))$ and $\dim(x(t^+))$, and $v$ is the impulse amplitude.

Note that by the definition of DK-STP, the dimension of the disturbance does not affect the nominal dimension of $x(t)$. If it is a state impulse, as long as the amplitude of $v$ is specified, say
 \begin{align}\label{3.4.3}
v=\mu\|(x(t^+)-x(t^-)\|_{{\cal V}},
\end{align}
 where $\mu>0$ is a constant parameter, the impact of the disturbance can be precisely estimated.

In the nonlinear case, the projected impulse
$$
\Pi^{\ell}_n\eta_k(t)
$$
can be used to estimate its impact on the nominal system, which has nominal dimension $n$.

\section{Differential Geometry on $\Omega$}

\subsection{Smooth Functions}

\begin{dfn}\label{d.30.1.1} $h(x):\Omega\ra \R$ is called a smooth ($C^{\infty}$) function, if $h(x)|_{\Omega^n}$, $n=1,2,\cdots$ are smooth ($C^{\infty}$) functions. The set of smooth functions is denoted by $C^{\infty}(\Omega)$.
\end{dfn}

To construct a smooth function on $\Omega$, we need some preparations.

\begin{lem}\label{l.30.1.2} \cite{che19} Assume $x\in \R^m$, $y\in \R^n$. $x\lra y$, if and only if, there exists a unique $z\in \R^s$, where $s=\gcd(m,n)$, such that $x=z\otimes \J_{m/s}$ and $y=z\otimes \J_{n/s}$.
\end{lem}

\begin{lem}\label{l.30.1.3}
\begin{align}\label{30.1.1}
\Pi^{km}_t(x\otimes \J_k)=\Pi^m_t(x),\quad x\in \R^m.
\end{align}
\end{lem}

\noindent{\it Proof.} Denote
$$
p=\lcm(t,km),\quad q=\lcm(t,m),\quad p=rq.
$$
Then we have
$$
\begin{array}{l}
\Pi^{km}_t(x\otimes \J_k)=\frac{t}{p}\left(I_t\otimes \J^{\mathrm{T}}_{p/t}\right)\left(I_{km}\otimes \J_{p/(km)}\right)(x\otimes \J_k)\\
 =\frac{t}{p}\left(I_t\otimes \J^{\mathrm{T}}_{p/t}\right)\left(I_k\otimes I_m\otimes \J_{p/(km)}\right)((x\otimes \J_k)\otimes 1)\\
 =\frac{t}{p}\left(I_t\otimes \J^{\mathrm{T}}_{p/t}\right)\left(x\otimes \J_k)\otimes \J_{p/km}\right)\\
 =\frac{t}{p}\left(I_t\otimes \J^{\mathrm{T}}_{p/t}\right)\left( x\otimes \J_{p/m}\right)\\
 =\frac{t}{p}\left( I_t\otimes \J^{\mathrm{T}}_{q/t}\otimes \J^{\mathrm{T}}_r\right)\left( x\otimes \J_{q/m}\otimes \J_{r}\right)\\
 =\frac{tr}{p}\left(I_t\otimes\J^{\mathrm{T}}_{q/t}\right)\left(x\otimes \J_{q/m}\right)\\
 =\frac{t}{q}\left(I_t\otimes \J^{\mathrm{T}}_{q/t}\right)\left(I_m\otimes\J_{q/m}\right)x
 =\Pi^m_t x.
 \end{array}
 $$
\hfill $\Box$

\begin{thm}\label{t.30.1.4} Let $x\in \R^m$, $y\in \R^n$, and $x\lra y$. Then for any $t\in \Z_+$
\begin{align}\label{30.1.2}
\Pi^{m}_t(x)=\Pi^n_t(y).
\end{align}
\end{thm}

\noindent{\it Proof.}
Let $s=\gcd(m,n)$. Using Lemma \ref{l.30.1.2}, there exists a $z\in \R^s$ such that
$$
x=z\otimes \J_{m/s} ,\quad y=z\otimes \J_{n/s}.
$$
Using Lemma \ref{l.30.1.3},
$$
\Pi^m_t(x)=\Pi^s_t(z)=\Pi^n_t(y).
$$
\hfill$\Box$

\begin{dfn}\label{d.30.1.5} Let $h(x)\in C^{\infty}(\R^n)$ be a smooth function. Extend $h(x)$ to $\R^{\infty}$ as follows:
\begin{align}\label{30.1.3}
H(y):=
\begin{cases}
h(y),\quad y\in \R^n,\\
h(\Pi^m_n(y)),\quad y\in \R^m,\quad m\neq n.
\end{cases}
\end{align}
\end{dfn}

Using Theorem \ref{t.30.1.4}, the following proposition is obvious.

\begin{prp} \label{p.30.1.6} Define
\begin{align}\label{30.1.4}
H(\bar{y}) :=H(y),\quad y\in \R^{\infty},
\end{align}
then $H$ is smooth function on $\Omega$, i.e., $H\in C^{\infty}(\Omega)$.
\end{prp}

\begin{rem}\label{r.30.1.7} To see $H$ is properly defined, we need to show that the $H(\bar{y})$ defined by (\ref{30.1.4}) is independent of the choice of $y\in \bar{y}$. This fact is ensured by Theorem \ref{t.30.1.4}.
\end{rem}

\subsection{Vector Fields}

Let $f(x)\in V(\R^n)$ be a smooth $C^{\infty}$ vector field. Consider a principal filter
${\cal F}_n:={\cal F}(\Omega^n)$. Define a vector field $F(x)$ as follows.
$$
F(y)=
\begin{cases}
f(y),\quad y\in \R^n,\\
\Pi^n_mf(\Pi^m_n y),\quad y\in {\cal F}_n.
\end{cases}
$$

\begin{prp}\label{p.30.2.1} Assume $n|m$, define
\begin{align}\label{30.2.1}
F(\bar{y})=F(y),\quad \bar{y}\in \Omega.
\end{align}
Then $F$ is a properly defined vector field. That is, if $y_1\lra y_2$, then $F(y_1)\lra F(y_2)$. Precisely, assume $y_1\in \R^{pn}$, $y_2\in \R^{qn}$, $s=\gcd(p,q)$, and $y_1\lra y_2$, then
\begin{align}\label{30.2.2}
F(y_1)\otimes \J_{p/s}=F(y_2)\otimes \J_{q/s}.
\end{align}
\end{prp}

\noindent{\it Proof.} According to Theorem \ref{t.30.1.4},
$$
z:=F(\Pi^{pn}_n(y_1))=F(\Pi^{qn}_n(y_2)).
$$
A straightforward computation shows that
$$
\Pi^n_{pn}z=z\otimes \J_p,\; \Pi^n_{qn}z=z\otimes \J_q.
$$
The conclusion follows.
\hfill $\Box$

\begin{rem}\label{r.30.2.2}
\begin{itemize}
\item[(i)] From the proof of the Proposition \ref{p.30.2.1} one sees that if $f(x)\in V(\R^n)$ (equivalently, $f(x)\in V(\Omega^n)$), then the vector field $F(\bar{x})\in V(\Omega)$ is only defined on $T(\Omega^{kn})$, $k=1,2,\cdots$. Moreover,
\begin{align}\label{30.2.3}
F(x\otimes \J_k)=f(x)\otimes \J_k,\quad k=1,2,\cdots.
\end{align}

\item[(ii)] For each $\bar{x}_0\in \Omega$, assume $\dim(\bar{x}_0)=s$ and $\lcm(s,n)=t$. Then $F(\bar{x}_0)$ is defined on
$T(\Omega^{kt})$, $k=1,2,\cdots$, which is a filter of $T(\bar{x}_0)=\bigcup_{k=1}^{\infty}\R^{ks}$.

\item[(iii)] Assume for each $F(x)\in V(\Omega)$, there is an $\Omega^{n_0}$ such that
\begin{align}\label{30.2.4}
F(x)|_{\Omega^{n_0}}\in V(\Omega^{n_0}),
\end{align}
i.e., it is defined at each $x\in \Omega^{n_0}$. Moreover, if $n_0$ is the smallest $n$ satisfying (\ref{30.2.4}), then $F(x)$ is said to be generated by $f(x)=F(x)|_{\Omega^{n_0}}$. In this paper, we are concerned only with such smooth vector fields.
 \end{itemize}
\end{rem}

 From the above argument, we have the following result.

\begin{prp}\label{p.30.2.3}
Assume $F(x)\in V(\Omega)$ is generated by $f(x)\in V(\Omega^n)$. $f_m(x)=F(x)|_{\Omega^m}$, where $n|m$.
Let $x_0\in \Omega^n$ and the integral curve of $f(x)$ with initial point $x(0)=x_0$ is $\Phi^{f(x)}_t(x_0)$, then the integral curve of
$f_m(y)$ with initial value $y_0=x_0\otimes \J_{m/n}$ is
\begin{align}\label{30.2.5}
\Phi^{f_m(y)}_t(y_0)=\Phi^{f(x)}_t(x_0)\otimes \J_{m/n}.
\end{align}
\end{prp}

\begin{rem}\label{r.30.2.3} Proposition \ref{p.30.2.3} is essential for our approach in the sequel. It means a vector field $f(x)$ on $\Omega^n$ can be lifted to $\Omega^{kn}$ as $F(y)$, and $F(x\otimes \J_k)=f(x)\otimes \J_k$. Moreover, the integral curve $\Phi^{f}_t(x_0)$ can be lifted to $\Phi^{F}_t(y_0)$ for $y_0=x_0\otimes \J_k$, such that
\begin{align}\label{30.2.6}
\Phi^{F}_t(x_0\otimes\J_k)=\Phi^f_t(x_0)\otimes \J_k.
\end{align}
In one word, the integral curve of $f$ is a by-product of the integral curve of its lifting vector field $F$. This observation makes the definition of vector field on $\Omega$ meaningful. Note that the vector fields on $\Omega$ are defined only over a filter of tangent bundles.
\end{rem}

\subsection{Distribution}

\begin{dfn}\label{d.30.3.1} Let $f_i(x)\in V(\Omega)$, $i\in [1,s]$ be a finite set of vector fields. Moreover, $f_i$ are generated by $f_i^0\in V(\Omega^{n_i}$.
Set $n=\lcm(n_1,\cdots,n_s)$. Then the distribution
\begin{align}\label{30.3.1}
D(x)=\Span\{f_i(x)\;|\;i\in[1,s]\},
\end{align}
which is defined on $\Omega^{kn}$ $k=1,2,\cdots$.
\end{dfn}

We refer to \cite{chepr} for constructing co-vector fields, co-distributions, and tensor fields over $\Omega$.

\section{Control Design of Dimension-Varying Systems}

\subsection{Dimension-Varying Control Systems}

\begin{dfn}\label{d4.1.1} Let $D:=\{1,2,\cdots,s\}\subset \Z_+$ be a finite set, and $\sigma:[0,\infty)\ra D$ be a right continuous piecewise constant function. A dimension-varying control system is
\begin{align}\label{4.1.1}
\begin{cases}
\dot(x)(t)=G_{\sigma(t)}(x(t),u)\subset \R^{n_{\sigma(t)}},\\
y(t)=h(\Pi^{n_{\sigma(t)}}_qx(t)),
\end{cases}
\end{align}
where $u=(u_1(t),\cdots,u_m(t))$ is the set of inputs (i.e., controls), $h:\R^q\ra \R^p$ is a (smooth) mapping, and $y(t)\in \R^p$ are the constant dimensional outputs.
\end{dfn}

As aforementioned, we use $\Omega$ as a common state space for any DV(C)S. But for a particular DV(C)S, only a sub-lattice of $\{\Omega^n\;|\;n=1,2,\cdots\}$ is required. We give an example to explain this.

\begin{exa}\label{e.4.1.2} Consider a DVCS, where $D=\{1,2,3\}$, and $G_1(x,u)\subset \R^2$, $G_2(x,u)\subset \R^3$, and $G_3(x,u)\subset \R^5$, then a sub-lattice as depicted in Figure \ref{Fig.4.1}, which is called a Hasse diagram. This sublattice is sufficient to capture the evolution of our switching system.

\begin{figure}
\begin{center}
\setlength{\unitlength}{15mm}
\begin{picture}(6,4)
\thicklines
\put(1,1){\circle*{0.2}}
\put(3,1){\circle*{0.2}}
\put(5,1){\circle*{0.2}}
\put(2,2){\circle*{0.2}}
\put(3,2){\circle*{0.2}}
\put(4,2){\circle*{0.2}}
\put(3,3){\circle*{0.2}}
\put(1,1){\line(1,1){2}}
\put(1,1){\line(2,1){2}}
\put(5,1){\line(-1,1){2}}
\put(5,1){\line(-2,1){2}}
\put(3,1){\line(1,1){1}}
\put(3,1){\line(-1,1){1}}
\put(3,2){\line(0,1){1}}
\put(0.9,0.6){$\Omega^{2}$}
\put(2.9,0.6){$\Omega^{3}$}
\put(4.9,0.6){$\Omega^{5}$}
\put(1.5,2){$\Omega^{6}$}
\put(4.2,2){$\Omega^{15}$}
\put(2.5,2){$\Omega^{10}$}
\put(2.8,3.2){$\Omega^{30}$}
\end{picture}\label{Fig.4.1}
\end{center}
\centerline{Fig.4.1~Sub-Lattice as a State Space}
\end{figure}

Consider $G_1(x,u)\subset \Omega^2$, It is expressed in $\Omega^6$ as $G_1(\Pi^6_2 y,u)\otimes \J_3$, $y\in\Omega^6$. And
it is expressed in $\Omega^{30}$ as $G_1(\Pi^{30}_2 z,u)\otimes \J_{15}$, $z\in\Omega^{30}$. etc.

\end{exa}

\begin{exa}\label{e.4.1.3}
\begin{itemize}
\item[(i)]
Consider a dimension-varying affine nonlinear control system
\begin{align}\label{4.1.2}
\begin{array}{ll}
\Sigma_i:&
\begin{cases}
\dot{x}(t)=f_i(x)+\dsum_{k=1}^mg^i_k(x)u_k,\\
y=h(\Pi^{n_i}_p x),\quad x\in \R^{n_i}.
\end{cases}\\
~&y\in \R^p,\quad i\in [1,s].
\end{array}
\end{align}

\item[(ii)]
Consider a dimension-varying linear control system
\begin{align}\label{4.1.3}
\begin{array}{ll}
\Sigma_i: &
\begin{cases}
\dot{x}(t)=A_i x(t)+\dsum_{k=1}^mb^i_k u_k,\\
y(t)=H\ttimes x(t), \quad x\in \R^{n_i},
\end{cases}\\
~~&i\in [1,s],
\end{array}
\end{align}
where $A_i\in {\cal M}_{n_i\times n_i}$, $b^i_j\in \R^{n_i}$, $H\in {\cal M}_{p\times q}$, $i\in [1,s]$.
\end{itemize}
\end{exa}

\subsection{Controllability}

\begin{dfn}\label{d.4.2.1} A dimension-varying control system is said to be completely controllable, if there exists a proper control $u(t)$, which can drive the trajectory from $x_0=x(0)\in \R^{n_i}$ to $x_e=x(T)$ ($T<\infty$) continuously.
\end{dfn}

\begin{rem}\label{r.4.2.2}
\begin{itemize}
\item[(i)] The continuity is under the distance induced topology ${\cal T}_d$ of $\R^{\infty}$. That is, in space $\Omega$. In the natural topology ${\cal T}_n$, this is impossible. See Figure \ref{Fig.4.2} for the continuous trajectory $\widetilde{AD}$ in $\Omega$, ignoring the segment $pq$.

\item[(ii)] Since discontinuous switches are practically with certain error, controllability seems questionable. If the controllability is so-called practically controllable, i.e., controlling the trajectory into a small neighborhood, the discontinuous trajectories become meaningful. The segment $pq$ in Figure \ref{Fig.4.2} shows a discontinuous switching.
\end{itemize}
\end{rem}

\begin{figure}
\begin{center}
\setlength{\unitlength}{10mm}
\begin{picture}(5,7)
\thinlines
\put(2,0){\line(1,1){3}}
\put(2,0){\line(0,1){2}}
\put(0,2){\line(0,1){2}}
\put(0,2){\line(1,0){2}}
\put(0,4){\line(1,1){3}}
\put(5,3){\line(0,1){2}}
\put(5,5){\line(-1,0){2}}
\put(5,5){\line(-1,-1){3}}
\put(3,5){\line(0,1){2}}
{\color{red}
\put(1.9,5){\vector(1,-3){0.45}}
}
\put(1.7,5){$p$}
\put(2.2,3.3){$q$}
\thicklines
\put(0,2){\line(1,1){3}}
\put(5,5){\line(-1,-1){3}}
\put(4,6){${\cal H}\subset\Omega$}
\put(0.5,4){$\Omega^m$}
\put(0.8,2.4){$\Omega^{m\wedge n}$}
\put(3,4.2){$\Omega^{n}$}
\put(3.5,2.5){$\Omega^s$}
\put(4.2,3.9){$\Omega^{n\wedge s}$}
\put(2.5,6.2){$A$}
\put(2.2,0.5){$D$}
\put(2.3,6.0){\circle*{0.2}}
\put(2.2,1){\circle*{0.2}}
\qbezier(2.3,6.0)(1.5,4.0)(2,4)
\put(1.5,3.8){$B$}
\qbezier(2,4)(3,3)(2.5,2.5)
\put(2.2,2.5){$C$}
\qbezier(2.5,2.5)(2.1,1.5)(2.2,1)
\end{picture}\label{Fig.4.2}
\end{center}

\centerline{Fig.4.2~~A Continuous Trajectory in $\Omega$ }
\end{figure}

Note that $\Omega$ is path-connected; precisely speaking, any two components $\Omega^m$ and $\Omega^n$ intersect at $\Omega^{m\wedge n}$, which is at least a one-dimensional subspace.

\begin{exa}\label{e.4.2.3} Consider $\Sigma=\{\Sigma_i\subset\R^{n_i}\;|\;i\in [1,s]\}$. Assume each mode $\Sigma_i$ is
$$
\dot{x}^i(t)=A_i+B_iu,\quad i\in [1,s],
$$
and $(A_i,B_i)$, $\forall i$ are completely controllable, then $\Sigma$ is completely controllable.

\end{exa}

In fact, the requirement that all the modes are controllable is too strong. We give a sufficient condition, which is obviously weaker than the above condition.

\begin{prp}\label{p.4.2.4} Consider a DVCS $\Sigma=\{\Sigma_i\subset \R^{n_i}\;|\;i\in [1,s]\}$. Consider the state space as the smallest sublattice generated by $\{\Omega^{n_i}\;|\; i\in[1,s]\}$, denoted by
\begin{align}\label{4.2.1}
{\cal H}={\cal L}(\Omega^{n_i}\;|\;i\in [1,s])
\end{align}
Assume $x_0\in \Omega^{n_p}\subset {\cal H}$ and $x_d\in \Omega^{n_q}\subset {\cal H}$. Then the state can be controlled from
$x_0$ to $x_d$, if
\begin{itemize}
\item[(i)] there exists a chain
$$
\Omega_1=\Omega^{n_p}\ra \Omega_2\ra \cdots \ra \Omega_r=\Omega^{n_q},
$$
such that $\Sigma_k$, $k\in [1,r-1]$ are partly controllable on
$$
{\Omega^k \bigcap \Omega^{k+1}}^{\perp},\quad k\in[1,r-1].
$$
\item[(ii)] $\Omega^r$ is completely controllable.
\end{itemize}
\end{prp}

\noindent{\it Proof.} The condition (i) ensures that the trajectory can go to ${\Omega_k \bigcap \Omega_{k+1}}$ sequentially.
The condition (ii) ensures the trajectory can go from a point on ${\Omega_{r-1} \bigcap \Omega_{r}}$ to $x_d$.
\hfill $\Box$

\begin{rem}\label{r.4.2.5} Observe Figure \ref{Fig.4.2}. The curve $\widetilde{AD}$ can be considered as a controlled trajectory of
$$
\Sigma=\{\Sigma_1\subset \Omega^{m}, \Sigma_2\subset \Omega^{n},\Sigma_3\subset \Omega^{s}\}.
$$
Assume $\Sigma_1$ is partly controllable on ${\Omega^{m\wedge n}}^{\perp}$, so it can drive the trajectory from $A$ to some point
$B\in \Omega^{m\wedge n}$. Similarly, $\Sigma_2$ can drive $B$ to certain point $C\in \Omega^{n\wedge s}$. Finally, $\Sigma_3$ is completely controllable, which drives $C$ to $D$.
\end{rem}

We give a numerical example for the controllability.

\begin{exa}\label{e.4.2.6} Consider $\Sigma=\{\Sigma_1\subset \Omega^2,\Sigma_2\subset \Omega^3\}$,
where
$$
\begin{array}{ll}
\Sigma_1:&\dot{x}(t)=\frac{1}{2}\begin{bmatrix}1&1\\-1&1\end{bmatrix}x(t)+\frac{1}{2}\begin{bmatrix}1\\-1\end{bmatrix}u(t);\\
\Sigma_2:&\dot{x}(t)=\begin{bmatrix}0&1&0\\0&0&1\\0&0&0\end{bmatrix}x(t)+\begin{bmatrix}0\\0\\1\end{bmatrix}u(t).\\
\end{array}
$$
Assume $x_0=(2,1)^{\mathrm{T}}\in \Omega^2$. Set
$y(t)=Tx(t)$, where
$$
T=\begin{bmatrix}
1&-1\\
1&1\\
\end{bmatrix}.
$$
Then we have
$$
\dot{y}(t)=\begin{bmatrix}1&1\\0&1\end{bmatrix}y(t)+\begin{bmatrix}1\\0\end{bmatrix}u(t).
$$
Choosing $u(t)=-(y_1+y_2)+v(t)$ yields
\begin{align}\label{4.2.2}
\begin{cases}
\dot{y}_1(t)=v,\\
\dot{y}_2(t)=y_2.
\end{cases}
\end{align}
Note that $\Omega^{2\wedge 3}=\Omega^1$, which intersects $\Omega^2$ at
$$
Omega^1=\{(r,r)^{\mathrm{T}}\;|\;r\in \R\}.
$$
Hence
$$y_1=(1,-1)\begin{bmatrix}x_1\\x_2\end{bmatrix}
$$
is the perpendicular direction of $\Omega^1$.
From (\ref{4.2.2}), it is clear that $y_1$ is controllable. Simply choose $v=-1$, then at $T=1$ we have
$$
y_1(T)=0;\quad y_2(T)=3e.
$$
That is
$$
x_1(T)=1.5e;\quad x_2=1.5e.
$$
That is, $x(T)\in \Omega^1$.

Secondly, since $\Sigma_2$ is completely controllable, a proper control can drive $x(T)\lra (1.5e,1.5e,1.5e)^{\mathrm{T}}$ to any point in $\Omega^3$.
\end{exa}

\subsection{Observability}

Consider system (\ref{4.1.2}) or (\ref{4.1.3}). We assume the observer $y$ is of constant dimension, the outputs are then considered as a mapping, consisting of a set of smooth functions over $\Omega$. Assume the function is generated by $h(x)\in C^{\infty}(\R^q)$, then it can be expressed as in (\ref{4.1.2}). In the linear case, it degenerates to a matrix, which is described by (\ref{4.1.3}).

The following result is an immediate consequence of the definition.

\begin{prp}\label{p.4.3.1}
\begin{itemize}
\item[(i)] Consider nonlinear control system \eqref{4.1.2}. Assume the output mapping $h$ is generated by $h_q\in C^{\infty}(\R^q)$. When $\Sigma_i$ is active, the corresponding output mapping becomes $h=h_q(\Pi^{n_i}_q x)$, where $x\in \R^{n_i}$.
\item[(ii)] Consider linear control system (\ref{4.1.3}). Assume the output mapping $h$ is generated by $H\in {\cal M}_{p\times q}$.
 When $\Sigma_i$ is active, the corresponding output mapping becomes $h=H\ltimes x=H\Psi_{q\times n_i} x$, where $x\in \R^{n_i}$.
\item[(iii)] If the switching signal is known and mode distinguishability is not considered, the system (\ref{4.1.2}) (or (\ref{4.1.3})) is completely observable provided that it is completely observable on all modes.
\end{itemize}
\end{prp}

\subsection{Stabilization}

Using the controllability results, the following sufficient condition holds.

\begin{prp}\label{p.4.4.1} Consider system (\ref{4.1.3}). Assume the minimum dwell time $\D_m>0$.
\begin{itemize}
\item[(i)] The initial points from $\Omega_0$ can be continuously stabilized to zero at finite time, if there exist a mode ($\Sigma_d\subset \Omega_e$), which is stabilizable to zero at arbitrary time period (say, $\Sigma_d$ is completely controllable),
and there exist a finite sequence
$$
\Omega_0\ra \Omega_1\ra\cdots\ra \Omega_k=\Omega_d,
$$
such that $\Sigma_i$, $i\in [0,s-1]$ are partly controllable over the subspaces
$$
\{\Omega_i\bigcap \Omega_{i+1}\}^{\perp},\quad i\in[0,s-1],
$$
and $\Omega_k$ is stabilizable to zero at arbitrary time.
\item[(ii)] If any $\Omega_i$, $i\in [1,s]$ satisfy the condition (i), the system is globally (continuously) stabilizable.
\end{itemize}
\end{prp}

Note that we require the trajectories to be continuous here. In the following section, the jump case is considered.

\section{Equivalent Switching System}

\subsection{Embedding System}

By embedding each mode into the common state space, a dimension-varying system can be converted into a classical switching system \cite{yan14}.

Denote by $\d(t)$ the Dirac function \cite{tay80}. Assume that at a switching moment $t$ we have a jump from $x(t^-)$ to $x(t^+)$. Denote
$$
\|x(t^+)-x(t^-)\|_{{\cal V}}=r\neq 0.
$$
Set the jump direction as a unit vector by
$$
\varphi(t)=\frac{1}{r}(x(t^+)-x(t^-)).
$$
Denote by $n(t)=\lcm(n_{t^-}, n_{t^+})$,
then the extended switching mode is defined by
\begin{align}\label{5.1.1}
\tilde{\Sigma}_k:=\d(t)\varphi(t)(x(t^+)-x(t^-))+\Sigma_k\subset \Omega^{n(t)}.
\end{align}

Denote $n=\lcm(n_i\;|\;i\in [1,s])$. Then $\tilde{\Sigma}_k$ can be merged into $\Omega^n$ by
\begin{align}\label{5.1.2}
\tilde{\Sigma}_k:=\left(\d(t)\varphi(t)(x(t^+)-x(t^-))+\Sigma_k\right)\otimes \J_{n/n(t)}.
\end{align}
In this way, a dimension-varying control system becomes a switched control system. The only difference is that each mode may have possible impulses.

Principally, the techniques for classical switching control systems remain available for the switching system in $\Omega^n\cong \R^n$ with modes described by (\ref{5.1.2}).

\subsection{Stabilization with Switching Jumps}

This subsection considers the stabilization again. Unlike the stabilization with continuous trajectories discussed in subsection 4.4, we now allow switching jumps. The Dirac function is used to formulate the jump. Moreover, if the additional switching impulses are considered as in \cite{xue22}, they can be formulated as additional Dirac functions. We use an example to describe it.

\begin{exa}\label{e.5.2.1}
Consider $\Sigma=\{\Sigma_1\subset \R^2,\Sigma_2\subset \R^3\}$, where
$$
\begin{array}{ll}
\Sigma_1:&
\dot{x}(t)=\begin{bmatrix}0&1\\0&0.1\end{bmatrix}x(t)+\begin{bmatrix}1\\0\end{bmatrix}u.\\
\Sigma_2:&
\dot{z}(t)=\begin{bmatrix}0.1&0&0\\0&0&1\\1&0&1\end{bmatrix}z(t)+\begin{bmatrix}0\\0\\1\end{bmatrix}u.
\end{array}
$$
It is obvious that neither $\Sigma_1$ nor $\Sigma_2$ is stabilizable.

We use the following state feedback to stabilize the controllable sub-states.
Consider $\Sigma_1$, we choose the control as
$u(t)=-x_1(t)-x_2(t)$. Consider $\Sigma_2$, we choose the control as $u(t)=-z_1(t)-z_2(t)-3z_3(t)$.
Then the closed-loop system becomes
$$
\begin{array}{ll}
\Sigma_1:&
\dot{x}(t)=\begin{bmatrix}0&-1\\0&0.1\end{bmatrix}x(t).\\
\Sigma_2:&
\dot{z}(t)=\begin{bmatrix}0.1&0&0\\0&0&1\\0&-1&-2\end{bmatrix}z(t).
\end{array}
$$
Assume the switching jumps are nearest jumps.
That is,
$$
\begin{array}{l}
\begin{bmatrix}
z_1\\
z_2\\
z_3\\
\end{bmatrix}(t^+)=\Pi^2_3
\begin{bmatrix}
x_1\\
x_2\\
\end{bmatrix}(t^-),\\
\begin{bmatrix}
x_1\\
x_2\\
\end{bmatrix}(t^+)=\Pi^3_2
\begin{bmatrix}
z_1\\
z_2\\
z_3
\end{bmatrix}(t^-),\\
\end{array}
$$
where
$$
\begin{array}{l}
\Pi^2_3=\frac{1}{2}\left(I_3\otimes \J_2^{\mathrm{T}}\right)\left(I_2\otimes \J_3\right),\\
\Pi^3_2=\frac{1}{3}\left(I_2\otimes \J_3^{\mathrm{T}}\right)\left(I_3\otimes \J_2\right).\\
\end{array}
$$

Assume $x(0)=x_0=(5,6)^{\mathrm{T}}$. We consider two switching rules.
\begin{itemize}
\item[(i)] Fixed period switching:
Assume the sequence of switching moments are
$$
\begin{array}{l}
t_1=1,~t_2=3,~t_3=4,~t_4=6, \cdots,\\~~t_{2k+1}=t_{2k}+1,~t_{2(k+1)}=t_{2k+1}+2,\cdots
\end{array}
$$
The simulation result is shown in Figure \ref{Fig.5.1} by the black line.
\item[(ii)] Random switching: We assume the minimum dwell time is $0.5$ and the maximum dwell time is $2$. The simulation result is shown in Figure \ref{Fig.5.1} by a red (and thin) line.
\end{itemize}
Note that for both cases, a rigorous proof of the stability can also be done by using Kharitonov's theorem \cite{min89}.

\end{exa}

\begin{figure}
\begin{center}
\setlength{\unitlength}{10mm}
\begin{picture}(8,7)(-1,-1)
\thicklines
\put(0,-1){\vector(0,1){7}}
\put(-1,0){\vector(1,0){7.5}}
\put(-0.5,-0.5){$0$}
\put(0.8,-0.5){$1$}
\put(1.8,-0.5){$2$}
\put(2.8,-0.5){$3$}
\put(3.8,-0.5){$4$}
\put(4.8,-0.5){$5$}
\put(5.8,-0.5){$6$}
\put(-0.5,1){$1$}
\put(-0.5,2){$2$}
\put(-0.5,3){$3$}
\put(-0.5,4){$4$}
\put(-0.5,5){$5$}
\put(6.2,-0.5){$t$}
\put(-0.5,6.3){$\|x(t)\|_{{\cal V}}$}
\qbezier(0,5.52)(1,4.7)(2,1.9)
\qbezier(2,1.9)(3,0.8)(4,0.2)
\qbezier(4,0.2)(5,0.05)(6,0.02)
\thinlines
{\color{red}
\qbezier(0,5.52)(1,4.8)(2,2.6)
\qbezier(2,2.6)(3,1.2)(4,1.1)
\qbezier(4,1.1)(5,0.6)(6,0.4)
}

\end{picture}\label{Fig.5.1}
\end{center}

\centerline{Fig.5.1~~Norm of Trajectories}
\end{figure}

\subsection{Disturbance Decoupling}

The system considered is
$$
\Sigma=\{\Sigma_k\subset \R^{n_k}\;|\;k\in [1,s]\},
$$
where $\Sigma$ is formulated as
\begin{align}\label{5.3.1}
\begin{cases}
\dot{x}(t)=f_{\sigma(t)}(t)+\dsum_{j=1}^mg^j_{\sigma(t)}(t)u_j(t)+\xi(t)v,\\
y(t)=h(\Pi^{d(x(t))}_px(t)),
\end{cases}
\end{align}
where $\xi(t)\in \R^n$, $v$ is a disturbance, $h:\R^n\ra \R^p$ is the output, and $n=\lcm(n_k\;|\; k\in[1,s])$.
The disturbance decoupling problem (DDP) is to design proper controls such that the disturbance does not affect the outputs.

Combining the system lifting argument in Section 4 and the well-known results from nonlinear control theory \cite{isi95}, the following result is obvious.

\begin{prp}\label{p.5.3.1} The DDP for system $\Sigma$, described in (\ref{5.3.1}), is solvable if and only if the DDP is solvable on each mode, which is described as
\begin{align}\label{5.3.2}
\begin{cases}
\dot{x}(t)=f_{k}(t)+\dsum_{j=1}^mg^j_{k}(t)u_j(t)+\Pi^n_{n_k}\xi(t)v,\\
y(t)=h(\Pi^{n_k}_px(t)),\quad k\in [1,s].
\end{cases}
\end{align}
\end{prp}

\noindent{\it Proof.} The only thing we have to claim is that, since $\xi(t)-\Pi^n_{n_k}\xi(t)$ is orthogonal to $\R^{n_k}$, it will not affect the states of $\Sigma_{k}$. Similarly, it does not affect the output $y(t)$.
\hfill $\Box$

We give an example.

\begin{exa}\label{e.5.3.2}
Consider $\Sigma=\{\Sigma_1\subset \R^2,\Sigma_2\subset \R^3\}$.
where $\Sigma$ is described as
$$
\begin{cases}
\dot{x}(t)=f_{\sigma(t)}(x(t))+g_{\sigma(t)}(x(t))u+\xi(t)v,\\
y=h(\Pi^{d(x(t))}_nx(t)),
\end{cases}
$$
where $n=\lcm(2,3)=6$, $\xi(t)\in V(\R^6)$, $v$ is the disturbance,
$$
f_1(x)=\begin{bmatrix}
x_2\\2x_1/3\end{bmatrix};\;g_1(x)=\begin{bmatrix}
-1-x_2\\1\end{bmatrix};
$$
$$
f_2(z)=\begin{bmatrix}
z_2-z_3\\z_1\\-z_1\end{bmatrix};\;g_2=\begin{bmatrix}
z_3-z_2\\1\\-1\end{bmatrix};
$$
$$
\begin{array}{l}
h(w)=w_1+w_2+w_3+w_4+w_5+w_6+w_1w_2,\\
\xi(w)=(1+w_5^2, -1-w^2_6,0,-1-w_2,1,w_1)^{\mathrm{T}},
\end{array}
$$
$\sigma(t):[0,\infty)\ra \{1,2\}$ is a right continuous piecewise constant switching function.

Restrict $h$ to subspace $\Omega^2$ , it becomes
$$
h_1(x)=x_1+x_2+\frac{1}{3}x_1^2.
$$
Then
$$
dh_1(x)=(1+\frac{2}{3}x_1,1).
$$
The largest $(f_1,g_1)$-invariant subspace, containing in $\ker(h_1)$ is (we refer to \cite{isi95} for calculating this),
$$
V_1=\Span\left\{
\begin{bmatrix}
-1\\
1+\frac{2}{3}x_1\\
\end{bmatrix}\right\}.
$$

Restrict $h$ to subspace $\Omega^3$ , it becomes
$$
h_2(z)=z_1+z_2+z_3+\frac{1}{2}z_1^2.
$$
Then
$$
dh_2(z)=(1+z_1,1,1).
$$
The largest $(f_2,g_2)$-invariant subspace, containing in $\ker(h_2)$ is
$$
V_2=\Span\left\{
\begin{bmatrix}
-1\\
1+z_1\\
0\\
\end{bmatrix},
\begin{bmatrix}
-1\\
0\\
1+z_1\\
\end{bmatrix}
\right\}.
$$

The restriction of $\xi$ to $\Omega^2$ is
$$
\xi_1=\Pi^6_2\xi=\frac{1}{3}\left.\begin{bmatrix}
\xi^1+\xi^2+\xi^3\\
\xi^4+\xi^5+\xi^6\\
\end{bmatrix}\right|_{\begin{array}{l}w_1=w_2=w_3=x_1,\\w_4=w_5=w_6=x_2\end{array}}=0,
$$
which means it is perpendicular to $\Omega^2$. Hence it does not affect $h_1$.

The restriction of $\xi$ to $\Omega^3$ is
$$
\begin{array}{l}
\xi_2=\Pi^6_3\xi=\frac{1}{2}\left.\begin{bmatrix}
\xi^1+\xi^2\\
\xi^3+\xi^4\\
\xi^5+\xi^6\\
\end{bmatrix}\right|_{\begin{array}{l}w_1=w_2=z_1,\\w_3=w_4=z_2,\\w_5=w_6=z_3\end{array}}\\
=\frac{1}{2}
\begin{bmatrix}
0\\
-1-z_1\\
1+z_1\\
\end{bmatrix}
\end{array}
$$
It is obvious that
$$
\xi_2\in V_2.
$$
We conclude that the DDP of system $\Sigma$ is solvable.

\end{exa}

\section{Hierarchical Networks}

Consider a hierarchical large-scale network, say, the internet. Since it could be changed frequently, switching the mode from time to time seems impractical. A practical way to model such systems could be to choose a specific level of intersection nodes as the switching level. Variations in lower-level nodes are considered disturbances.

We describe this using a three-level network, as depicted in Figure \ref{Fig.6.1}. We choose the middle level as the switching modes. That is, the system $\Sigma$ is considered as a switching system
$$
\dot{x}(t)=f_{\sigma(t)}(x(t),u(t)),
$$
where $\sigma:t\ra [1,s]$ is a right continuous piecewise constant mapping.

The nominal model for $\Sigma^i$ is
\begin{align}\label{6.1}
\dot{x}^i(t)=f^i_{\sigma(t)}(x^i(t),u^i(t))\subset \R^{n_i},\quad i\in [1,s],
\end{align}
where the nominal dimension mode $f^i_{\sigma}$ could be the approximate mode for all $\Sigma^i_j$, $j\in [1,s]$, and $n_i$ could be the average of $n^i_j$, which is the dimension of $\Sigma^i_j$.

Now if $\Sigma^i_j$ is active, the approximate mode becomes
\begin{align}\label{6.2}
\begin{array}{l}
\dot{x}^i(t)=f^i_{\sigma(t)}(\Pi^{n^{i_j}}_{n^i}x^i_j(0)),u^i(t))\subset \R^{n_i},\\
x^i_j(t)=\Pi^{n_i}_{n^i_j}x^i(t),\quad i\in [1,s],
\end{array}
\end{align}

\begin{figure}
\begin{center}
\setlength{\unitlength}{6mm}
\begin{picture}(14,7)
\thicklines
\put(1,1){\line(0,1){2}}
\put(2,1){\line(0,1){2}}
\put(4,1){\line(0,1){2}}
\put(5,1){\line(0,1){2}}
\put(6,1){\line(0,1){2}}
\put(8,1){\line(0,1){2}}
\put(10,1){\line(0,1){2}}
\put(11,1){\line(0,1){2}}
\put(13,1){\line(0,1){2}}
\put(1,3){\line(1,0){3}}
\put(5,3){\line(1,0){3}}
\put(10,3){\line(1,0){3}}
\put(2.5,3){\line(0,1){2}}
\put(6.5,3){\line(0,1){3}}
\put(11.5,3){\line(0,1){2}}
\put(2.5,5){\line(1,0){9}}
\put(0.6,0.2){$\Sigma^1_1$}
\put(1.6,0.2){$\Sigma^1_2$}
\put(3.6,0.2){$\Sigma^1_{k_1}$}
\put(4.7,0.2){$\Sigma^2_1$}
\put(5.7,0.2){$\Sigma^2_2$}
\put(7.7,0.2){$\Sigma^2_{k_2}$}
\put(9.8,0.2){$\Sigma^s_1$}
\put(10.8,0.2){$\Sigma^s_2$}
\put(12.8,0.2){$\Sigma^s_{k_s}$}
\put(2.5,2){$\cdots$}
\put(6.5,2){$\cdots$}
\put(11.5,2){$\cdots$}
\put(1.7,4){$\Sigma^1$}
\put(5.7,4){$\Sigma^2$}
\put(10.7,4){$\Sigma^s$}
\put(6,6.2){{\large$\Sigma$}}
\put(8.5,4){$\cdots$}
\end{picture}\label{Fig.6.1}
\end{center}

\centerline{Fig.5.2~~A Hierarchical Network}
\end{figure}

To see that the dynamics of aggregated systems $\Omega^i_j$ can be approximated by their nominal system $\Omega^i$, we use a numerical example to demonstrate this. First, we introduce a linear approximation, which was introduced by \cite{che19b}, (see also \cite{che23} page 187-188).

\begin{prp}\label{p.6.1}
\begin{itemize}
\item[(i)] Consider a linear system
\begin{align}\label{6.3}
\dot{x}(t)=Ax(t),\quad x(0)=x_0,\;x(t)\in \R^n.
\end{align}
The least-squares approximated system in $\R^m$ is
\begin{align}\label{6.4}
\dot{z}(t)=A_{\pi}z(t),\quad z(0)=\Pi^n_m(x_0),\;z(t)\in \R^m,
\end{align}
where
\begin{align}\label{6.401}
A_{\pi}=
\begin{cases}
\Pi^n_m A (\Pi^n_m)^{\mathrm{T}}(\Pi^n_m(\Pi^n_m)^{\mathrm{T}})^{-1},\quad n\geq m,\\
\Pi^n_m A ((\Pi^n_m)^{\mathrm{T}}\Pi^n_m)^{-1}(\Pi^n_m)^{\mathrm{T}},\quad n< m.\\
\end{cases}
\end{align}
\item[(ii)] Consider a linear control system
\begin{align}\label{6.5}
\begin{cases}
\dot{x}(t)=Ax(t)+Bu,\quad x(t)\in \R^n,\; u\in \R^m,\\
y(t)=Cx(t),\quad y(t)\in \R^p.
\end{cases}
\end{align}
The least-squares approximated system in $\R^m$ is
\begin{align}\label{6.6}
\begin{cases}
\dot{z}(t)=A_{\pi}z(t)+\Pi^n_mBu,\quad z(t)\in \R^m,\\
y(t)=C_{\pi}z(t),
\end{cases}
\end{align}
where $A_{\pi}$ is as in (\ref{6.401}), and
\begin{align}\label{6.7}
C_{\pi}=
\begin{cases}
C(\Pi^n_m)^{\mathrm{T}} (\Pi^n_m(\Pi^n_m)^{\mathrm{T}})^{-1},\quad n\geq m,\\
C((\Pi^n_m)^{\mathrm{T}}\Pi^n_m)^{-1}(\Pi^n_m)^{\mathrm{T}},\quad n< m.\\
\end{cases}
\end{align}
\end{itemize}
\end{prp}

Now we are ready to present an example to show that the nominal system approximation of a set of near-dimensional systems is meaningful.

\begin{exa}\label{e.6.2} Consider a linear system
\begin{align}\label{6.8}
\dot{x}(t)=Ax(t),\quad x(0)=x_0,\; x(t)\in \R^n.
\end{align}
Its trajectory is $x(t)=e^{At}x_0$.
\begin{itemize}
\item[(i)] Assume we use $z(t)\in \R^{n-1}$ to approximate the trajectory of (\ref{6.8}). Using (\ref{6.401}), we set
$$
B=A_{\pi}=\Pi^n_{n-1} A (\Pi^n_{n-1})^{\mathrm{T}}(\Pi^n_{n-1}(\Pi^n_{n-1})^{\mathrm{T}})^{-1}.
$$
Set $x_0=500\J_{n}$, we have $z_0=\Pi^n_{n-1}x_0$. Then we calculate
$$
\begin{array}{l}
z(t)=e^{Bt}z_0,\\
\tilde{x}(t)=\Pi^{n-1}_nz(t).
\end{array}
$$
Finally, the relative error is defined as
$$
E(t)=\frac{\|\tilde{x}(t)-x(t)\|_{{\cal V}}}{\|x(t)\|_{{\cal V}}}.
$$
We choose $A$ so that it evolves slowly enough to simulate a stable system.
$$
\begin{array}{l}
A_1=0.001*I_n;\\
A_2=-0.001*diag(1,2,\cdots,n);\\
A_3=0.001*rand(n);\\
\end{array}
$$
Then the relative error $E(t)$ is shown in Table \ref{tb.6.1}.

\begin{table}[!htb]
\centering
\caption{$E(t)$ for Lower Dimension Approximation \label{tb.6.1}}
\vskip 2mm
\doublerulesep 0.5pt
\begin{tabular}{|c||c|c|c|c|}
\hline
$A\backslash t$&1&2&3&4\\
\hline
\hline
$A_1\backslash n-1(e-15)$&0.6799&0.6767&0.6783&0.6710\\
\hline
$A_1\backslash n-3(e-15)$&0.5611&0.5518&0.5566&0.5518\\
\hline
$A_1\backslash n-5(e-15)$&0.1391&0.1444&0.1550&0.1502\\
\hline
$A_2\backslash n-1$&0.0000&0.000&0.000&0.0000\\
\hline
$A_2\backslash n-3$&0.0000&0.0000&0.001&0.0001\\
\hline
$A_2\backslash n-5$&0.0000&0.0001&0.0001&0.0001\\
\hline
$A_3\backslash n-1$&0.0013&0.0026&0.0038&0.0049\\
\hline
$A_3\backslash n-3$&0.0013&0.0025&0.0037&0.048\\
\hline
$A_3\backslash n-5$&0.0012&0.0024&0.0036&0.0048\\
\hline
\hline
$A\backslash t$&$\cdots$&98&99&100\\
\hline
$A_1\backslash n-1(e-15)$&$\cdots$&0.6683&0.6771&0.6739\\
\hline
$A_1\backslash n-3(e-15)$&$\cdots$&0.5512&0.5484&0.5616\\
\hline
$A_1\backslash n-5(e-15)$&$\cdots$&0.1166&0.1413&0.1164\\
\hline
$A_2\backslash n-1$&$\cdots$&0.0007&0.0007&0.0007\\
\hline
$A_2\backslash n-3$&$\cdots$&0.0019&0.0019&0.0019\\
\hline
$A_2\backslash n-5$&$\cdots$&0031&0.0031&0.0032\\
\hline
$A_3\backslash n-1$&$\cdots$&0.0269&0.0269&0.0270\\
\hline
$A_3\backslash n-3$&$\cdots$&0.0262 &0.0262&0.0262\\
\hline
$A_3\backslash n-5$&$\cdots$&0.0250&0.0251&0.0251\\
\hline
\end{tabular}
\end{table}

\item[(ii)] Assume we use $z(t)\in \R^{n+1}$ to approximate the trajectory of (\ref{6.8}). Using (\ref{6.401}), we set
$$
C=A_{\pi}=\Pi^n_{n+1} A ((\Pi^n_{n+1})^{\mathrm{T}}\Pi^n_{n+1})^{-1}(\Pi^n_{n+1})^{\mathrm{T}}.
$$
Set $x_0=500\J_{n}$, we have $z_0=\Pi^n_{n+1}x_0$. Then we calculate
$$
\begin{array}{l}
z(t)=e^{Ct}z_0,\\
\tilde{x}(t)=\Pi^{n+1}_nz(t).
\end{array}
$$
Finally, the relative error $E(t)$ is shown in Table \ref{tb.6.2}.

\begin{table}[!htb]
\centering
\caption{$E(t)$ for Higher Dimension Approximation\label{tb.6.2}}\footnote{In this table $n+k$ or $n-k$ stands for the dimension of approximated system, $e-13$ etc. means the error is $a{10}^{-13}$}
\vskip 2mm
\doublerulesep 0.5pt
\begin{tabular}{|c||c|c|c|c|}
\hline
$A\backslash t$&1&2&3&$\cdots$\\
\hline
\hline
$A_1\backslash n+1(e-13)$&0.0086&0.0082&0.0091&0.1651\\
\hline
$A_1\backslash n+3(e-13)$&0.0022&0.0037&0.0052&0.0067\\
\hline
$A_1\backslash n+5(e-13)$&0.0036&0.0044&0.0058&0.0077\\
\hline
$A_2\backslash n+1$&0.0000&0.000&0.0000&0.0000\\
\hline
$A_2\backslash n+3$&0.0000&0.0000&0.0000&0.0001\\
\hline
$A_2\backslash n+5$&0.0000&0.0001&0.0001&0.0001\\
\hline
$A_3\backslash n+1$&0.0013&0.0025&0.0037&0.0048\\
\hline
$A_3\backslash n+3$&0.0012&0.0024&0.0035&0.0045\\
\hline
$A_3\backslash n+5$&0.0011&0.0022&0.0032&0.0042\\
\hline
\hline
$A\backslash t$&$\cdots$&98&99&100\\
\hline
\hline
$A_1/\backslash n+1(e-13)$&$\cdots$&0.16666&0.1692&0.1798\\
\hline
$A_1/\backslash n+3(e-13)$&$\cdots$&0.1616&0.1634&0.1657\\
\hline
$A_1/\backslash n+5(e-13)$&$\cdots$&0.1578&0.1601&0.1623\\
\hline
$A_2\backslash n+1$&$\cdots$&0.0007&0.0007&0.0007\\
\hline
$A_2\backslash n+3$&$\cdots$&0.0018&0.0018&0.0018\\
\hline
$A_2\backslash n+5$&$\cdots$&0.0028&0.0028&0.0029\\
\hline
$A_3\backslash n+1$&$\cdots$&0.0264&0.0264&0.0264\\
\hline
$A_3\backslash n+3$&$\cdots$&0.0246&0.0247&0.0247\\
\hline
$A_3\backslash n+5$&$\cdots$&0.0227&0.0227&0.0227\\
\hline
\end{tabular}
\end{table}

\end{itemize}
\end{exa}

\begin{rem}\label{r.6.3}
\begin{itemize}
\item[(i)] For linear system $\Sigma\subset \R^n$ the least-squares approximated system $\tilde{\Sigma}\subset \in \R^m$ means the trajectory of $\tilde{\Sigma}$ is the least-squares approximation to the trajectory of $\Sigma$.
\item[(ii)] For a nonlinear system $f\in V(\R^n)$, the ``error" of approximated $\tilde{f}\in V(\R^m)$ may be measured by
$$
\|\tilde{f}-f\|_{{\cal V}}.
$$
\end{itemize}
\end{rem}

\section{Illustrative Power-System Example: Generator Removal and Reconnection}
\label{sec:illustrative_power_system}

This section provides an illustrative power-system example. Its purpose is not to propose a new power-system stability method, but to show how the quotient-space construction developed above can be instantiated in a nontrivial engineering system in which the active dynamic dimension, the algebraic network constraints, and the segment-local operating point change simultaneously. All cross-dimensional comparisons below are performed on translated reduced representatives; raw DAE states are not compared directly in $d_{\mathcal V}$.

\subsection{System Setup and Switching Scenario}
\label{subsec:ps_setup_illustrative}

\begin{figure}[htp]
 \centering
 \includegraphics[width=0.7\columnwidth]{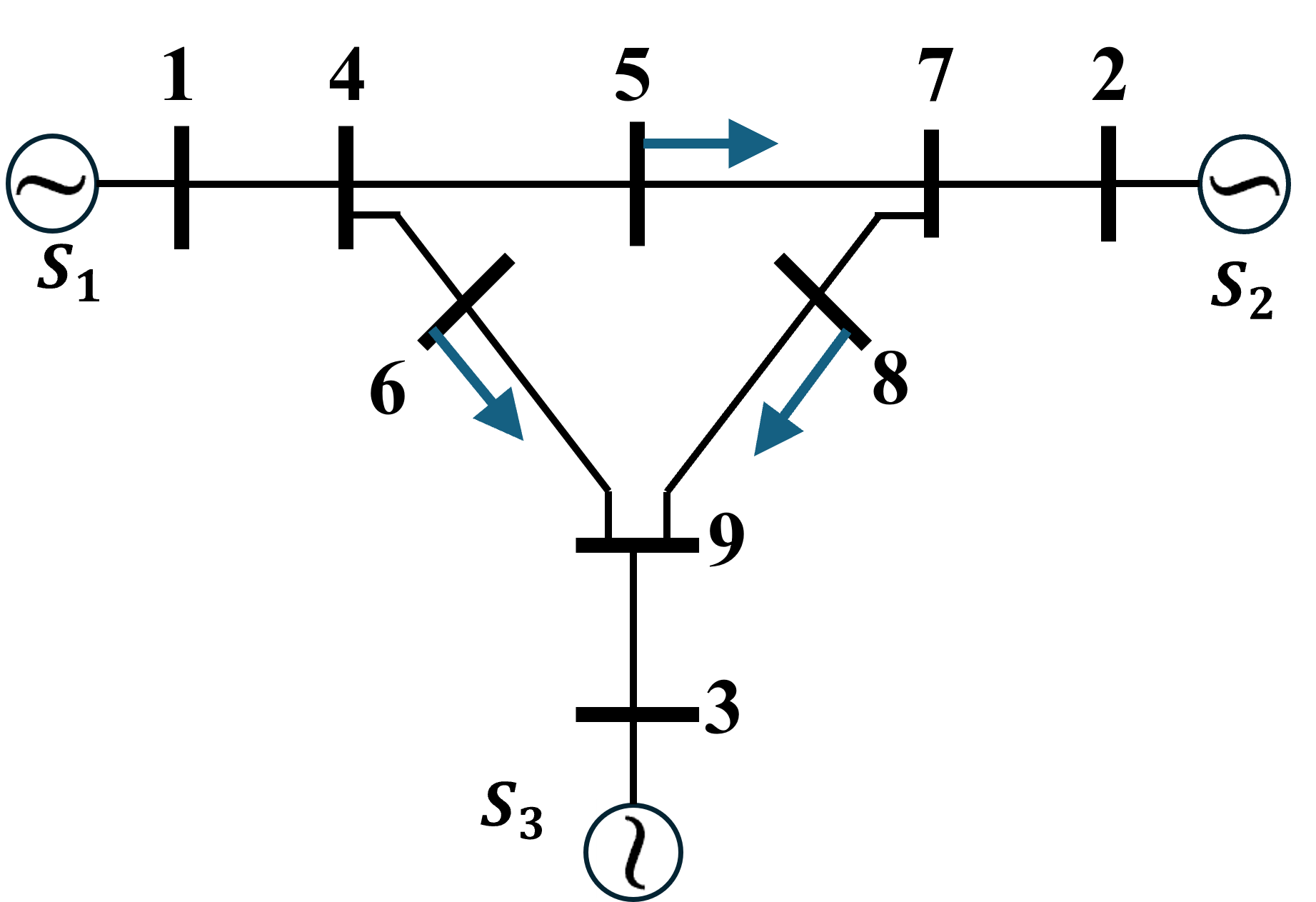}
 \caption{Illustrative three-machine nine-bus switching example. The finite schedule includes a load step at $t=1\,{\rm s}$, removal of $s_2$ at $t=2\,{\rm s}$, and reconnection of $s_2$ at $t=9\,{\rm s}$.}
 \label{fig:ill_one_line_diagram}
\end{figure}

Consider the classical three-machine nine-bus benchmark shown in Fig.~\ref{fig:ill_one_line_diagram}. Let $\mathcal S:=\{S_{123},S_{13}\}$ denote the set of active-machine modes, where
\begin{equation*}
S_{123}:=\{s_1,s_2,s_3\},\qquad S_{13}:=\{s_1,s_3\}.
\end{equation*}
For $S\in\mathcal S$, the superscript $(S)$ denotes the corresponding mode-wise reduced representative. The finite event schedule used in this example is $S_{123}\to S_{13}\to S_{123}$: a load step at $t=1.0\,{\rm s}$, removal of $s_2$ at $t=2.0\,{\rm s}$, and reconnection of $s_2$ at $t=9.0\,{\rm s}$. Hence the native dynamic dimension changes as
\begin{equation*}
6\;\longrightarrow\;4\;\longrightarrow\;6,
\end{equation*}
whereas, after removal of the common angle and common frequency freedoms, the internal representative dimension changes as
\begin{equation*}
4\;\longrightarrow\;2\;\longrightarrow\;4.
\end{equation*}

For each $S\in\mathcal S$, the original model is a mode-dependent power-system DAE. For quotient-space comparison, each mode is first locally reduced to a network-reduced ODE representative around a regular AC operating point. The detailed DAE equations, Kron-reduced network representation, and algebraic-elimination assumptions are collected in \ref{app:power_model_details}. The main text uses only the resulting reduced representatives and their translated forms.

\subsection{Construction of Translated Quotient Representatives}
\label{subsec:translated_representatives_illustrative}

\begin{figure}[htp]
 \centering
 \includegraphics[width=0.88\columnwidth]{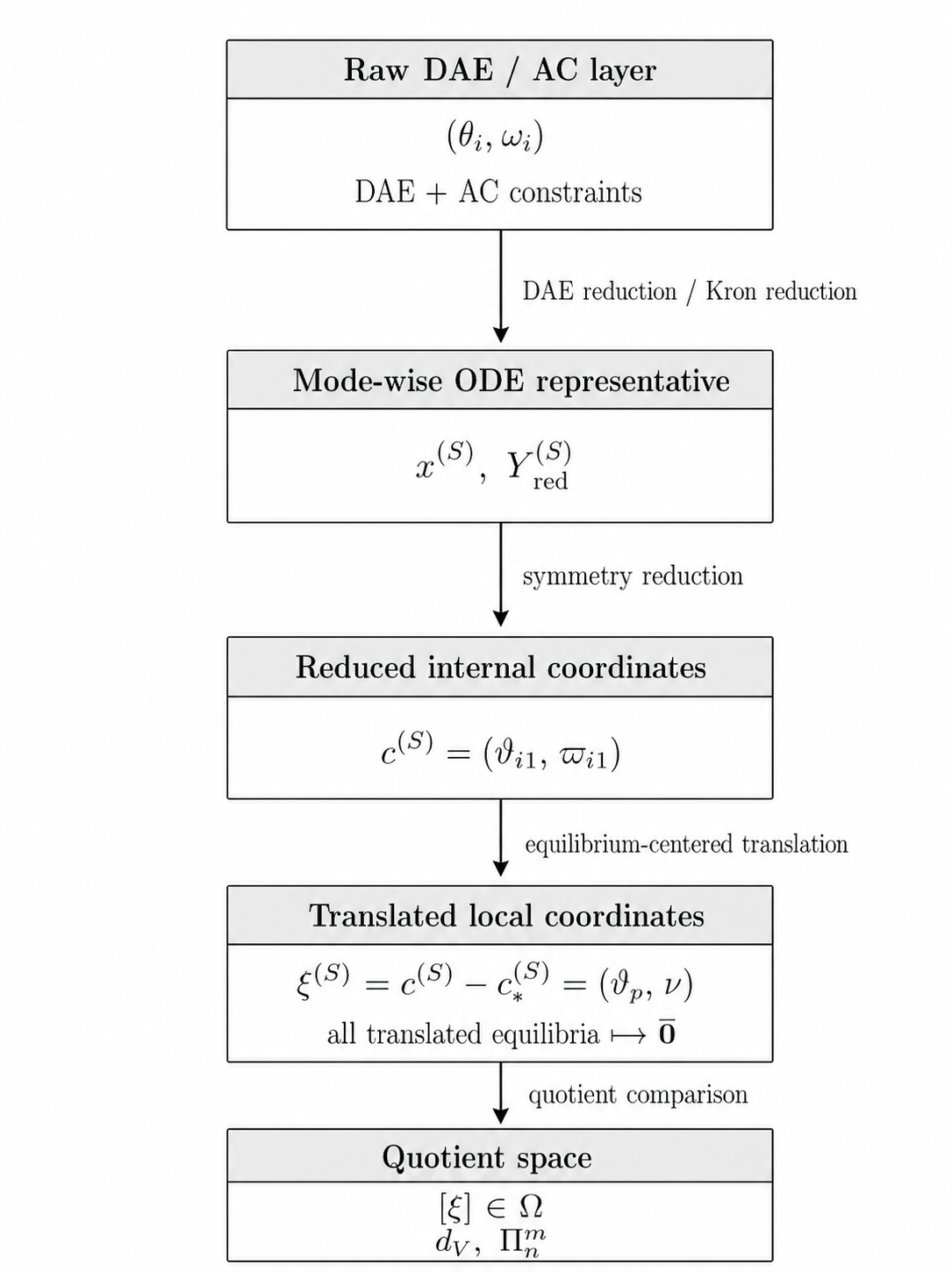}
 \caption{Construction of quotient-space representatives for the illustrative power-system example. The metric $d_{\mathcal V}$ is applied only after symmetry reduction and equilibrium-centered translation.}
 \label{fig:ill_coordinate_hierarchy}
\end{figure}

The construction follows the hierarchy in Fig.~\ref{fig:ill_coordinate_hierarchy}. Starting from raw DAE variables, one first obtains a mode-wise ODE representative, then removes the common rotational and frequency freedoms, and finally translates each reduced representative by its segment-local equilibrium. This is the point at which the power-system variables enter the quotient-space framework.

For a mode $S$, let $c^{(S)}$ denote the reduced internal representative and let $c_*^{(S)}$ denote its segment-local equilibrium. The translated representative is
\begin{equation*}
 \xi^{(S)}:=c^{(S)}-c_*^{(S)}.
\end{equation*}
For the three-machine and two-machine modes used here,
\begin{equation*}
 \xi^{(S_{123})}\in\mathbb R^4,
 \qquad
 \xi^{(S_{13})}\in\mathbb R^2.
\end{equation*}
Thus, each fixed-dimensional component retains its ordinary Euclidean reduced coordinates inside $\Omega^4$ or $\Omega^2$; the quotient geometry is invoked only when the representatives are compared across dimension-changing events. After this translation, all segment-local equilibria are represented by the same quotient object,
\begin{equation*}
 [0_4]=[0_2]=\bar 0\in\Omega.
\end{equation*}
This statement is a coordinate-level representation of the local equilibrium deviations. It is not a claim that the raw operating points of structurally different power-system modes are physically identical. The raw operating points, the reduced equilibrium representatives, and the event-state provenance are given in \ref{app:reduced_coordinates_details}.

Table~\ref{tab:ill_theory_case_map} summarizes how the main theoretical objects introduced earlier are instantiated in this example.

\begin{table}[!t]
 \centering
 \caption{Theory-to-example map for the illustrative power-system example.}
 \label{tab:ill_theory_case_map}
 \scriptsize
 \begingroup
 \setlength{\tabcolsep}{1.4pt}
 \renewcommand{\arraystretch}{1.03}
 \begin{tabular}{>{\raggedright\arraybackslash}p{0.29\columnwidth}
 >{\raggedright\arraybackslash}p{0.34\columnwidth}
 >{\raggedright\arraybackslash}p{0.29\columnwidth}}
 \toprule
 Earlier theoretical object/result & Power-system instantiation & Illustrated role \\
 \midrule
 Active-machine modes & $\mathcal S=\{S_{123},S_{13}\}$ & Structural modes in the finite schedule \\
 $\Omega^n\cong\mathbb R^n$ & $\Omega^4$ for $S_{123}$ and $\Omega^2$ for $S_{13}$ & Fixed-mode coordinates remain Euclidean \\
 Quotient equivalence & $[0_4]=[0_2]=\bar0$ & Common translated equilibrium object \\
 Finite sublattice & $\Omega^2\prec\Omega^4$ & Generated reduced dimensions only \\
 Projection/lift $\Pi_m^n$ & $\Pi_2^4$ and $\Pi_4^2$ & Removal/reconnection benchmarks \\
 DK-STP bridge interpretation & Same matrices in \eqref{eq:ill_projection_lift} & Bridge maps link multiplication and projection \\
 Switching gap $d_{\mathcal V}$ & Event-wise gaps in Table~\ref{tab:ill_event_compare} & Cross-dimensional discontinuity \\
 Lipschitz switching & $W_{\rm rem}^{\rm red}$ and $A_{\rm rec}$ & Finite gains; offsets explain displacement \\
 Dwell-time mechanism & $\tau_{13}=7\,{\rm s}$ & Conservative residence after removal \\
 \bottomrule
 \end{tabular}
 \endgroup
\end{table}

\subsection{Projection/Lift Benchmarks and Event-wise $d_{\mathcal V}$-Gaps}
\label{subsec:eventwise_metric_illustrative}

For this reduced benchmark, the relevant finite sublattice is
\begin{equation*}
 \mathcal L_{\rm red}=\{\Omega^2,\Omega^4\},
 \qquad
 {\rm lcm}(4,2)=4.
\end{equation*}
Since $2\mid4$, the two-machine reduced component is the lower element $\Omega^2\prec\Omega^4$ of the finite generated sublattice used by this example. This gives a concrete instance of the lattice structure introduced above.
The projection from the three-machine reduced space to the two-machine reduced space and the lift in the reverse direction are
\begin{equation}
 \Pi_2^4=\frac12
 \begin{bmatrix}1&1&0&0\\0&0&1&1\end{bmatrix},
 \qquad
 \Pi_4^2=
 \begin{bmatrix}1&0\\1&0\\0&1\\0&1\end{bmatrix}.
 \label{eq:ill_projection_lift}
\end{equation}
Thus, for $x\in\mathbb R^4$ and $z\in\mathbb R^2$, the metric evaluation reduces to
\begin{equation}
 d_{\mathcal V}(x,z)=\frac12\left\|x-z\otimes\mathbf 1_2\right\|_2.
 \label{eq:ill_dv_rule}
\end{equation}
The matrices in \eqref{eq:ill_projection_lift} also realize the bridge-matrix role that appears in the DK-STP interpretation of cross-dimensional matrix multiplication. In this example, they are used as quotient-geometric benchmarks, not as a claim that the physical tripping or reconnection policy must follow the nearest projection rule.

Applying \eqref{eq:ill_dv_rule} to the translated pre- and post-event representatives gives the event-wise comparison in Table~\ref{tab:ill_event_compare}. The realized event maps are affine; their linear parts have finite $d_{\mathcal V}$-induced gains, while the affine offsets account for the observed post-event displacement. The detailed vectors, induced gains, and arithmetic are reported in \ref{app:event_metric_recalculation}.

\begin{table}[!t]
 \centering
 \caption{Event-wise $d_{\mathcal V}$-gap comparison in translated reduced coordinates.}
 \label{tab:ill_event_compare}
 \scriptsize
 \begingroup
 \setlength{\tabcolsep}{2.6pt}
 \begin{tabular}{@{}lccc@{}}
 \toprule
 Event & $\delta^{\rm shift}$ & Benchmark & Residual \\
 \midrule
 Removal & $4.0056\times10^{-3}$ & $7.83\times10^{-5}$ & $4.0049\times10^{-3}$ \\
 Reconnection & $7.8869\times10^{-3}$ & $0$ & $7.8869\times10^{-3}$ \\
 \bottomrule
 \end{tabular}
 \endgroup
\end{table}

The table illustrates the main role of the quotient metric in this example: it separates the intrinsic projection/lift cost from the mismatch introduced by the realized switching policy and the change of segment-local equilibrium. During the removal event, the intrinsic projection cost is small relative to the realized shifted gap. In the reconnection event, the exact lift has zero intrinsic quotient cost, so the observed gap is entirely induced by the reinsertion policy and the affine offset. The zero benchmark cost in the reconnection event is a concrete instance of the quotient equivalence: $z\in\mathbb R^2$ and $z\otimes\mathbf 1_2\in\mathbb R^4$ represent the same element of $\Omega$. Equivalently, the residual column can be read as the quotient-metric magnitude of the event-induced deviation relative to the geometric benchmark.

\subsection{Dwell-Time Consistency with the Stability Mechanism}
\label{subsec:dwell_time_consistency_illustrative}

The sufficient-dwell-time argument developed above suggests that, when mode-wise Lyapunov decrease is available and switching gaps are uniformly bounded, stability can be ensured if the minimum dwell time is sufficiently large. The present example is not intended to compute the sharp dwell-time threshold. Instead, it illustrates the same mechanism on a realized finite switching schedule.

After the removal of $s_2$ at $t=2\,{\rm s}$, the system remains in the two-machine mode $S_{13}$ until $s_2$ is reconnected at $t=9\,{\rm s}$. Thus the residence time of the reduced mode is
\begin{equation*}
 \tau_{13}=9-2=7\,{\rm s}.
\end{equation*}
This value is the residence time of the dimension-reduced mode $S_{13}$, not a computed minimum dwell-time bound for the entire hybrid schedule. It is deliberately chosen to be conservative. Since the load step does not change the active dynamic dimension, the dwell-time interpretation here focuses on the dimension-changing removal--reconnection interval. The interval gives the active mode time to dissipate the transient energy introduced by the removal event before the next structural switching occurs. In this sense, the selected schedule is consistent with the sufficient-dwell-time mechanism of the general theory.

\begin{figure}[!t]
 \centering
 \includegraphics[width=0.78\columnwidth]{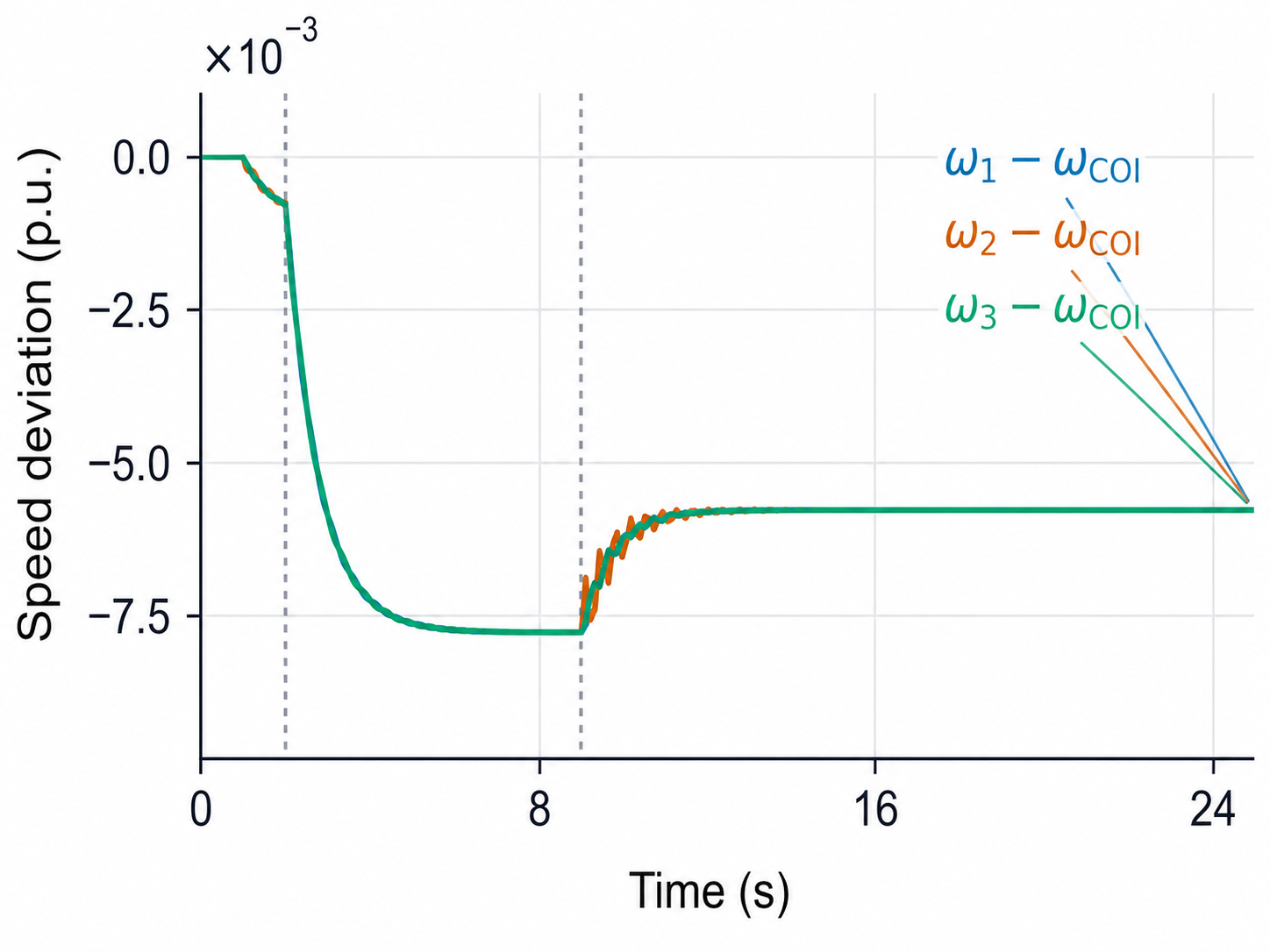}\\[-2pt]
 {\scriptsize (a) COI-relative speed deviations}\\[2pt]
 \includegraphics[width=0.78\columnwidth]{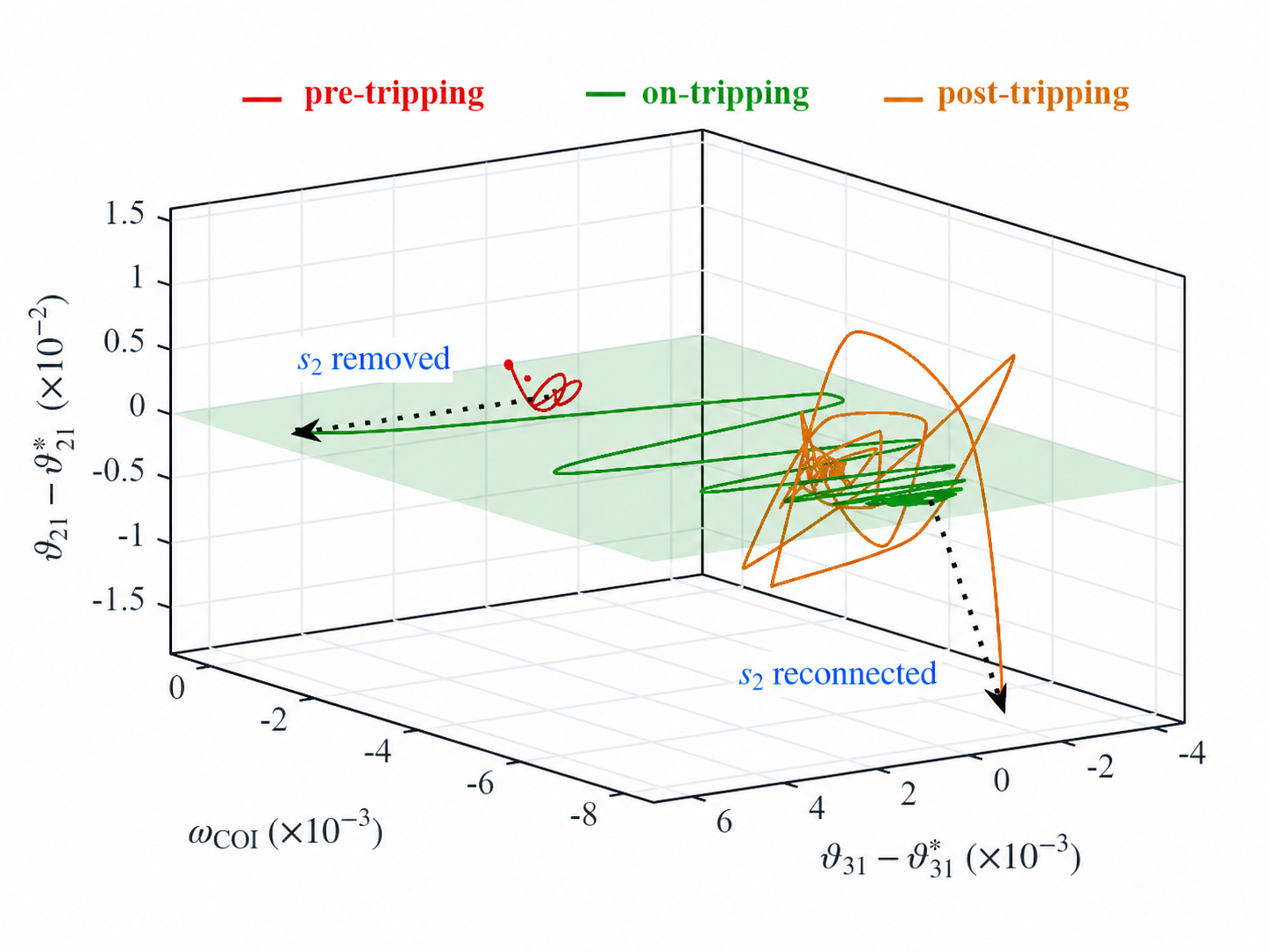}\\[-2pt]
 {\scriptsize (b) Reduced switching trajectory}\\[2pt]
 \includegraphics[width=0.78\columnwidth]{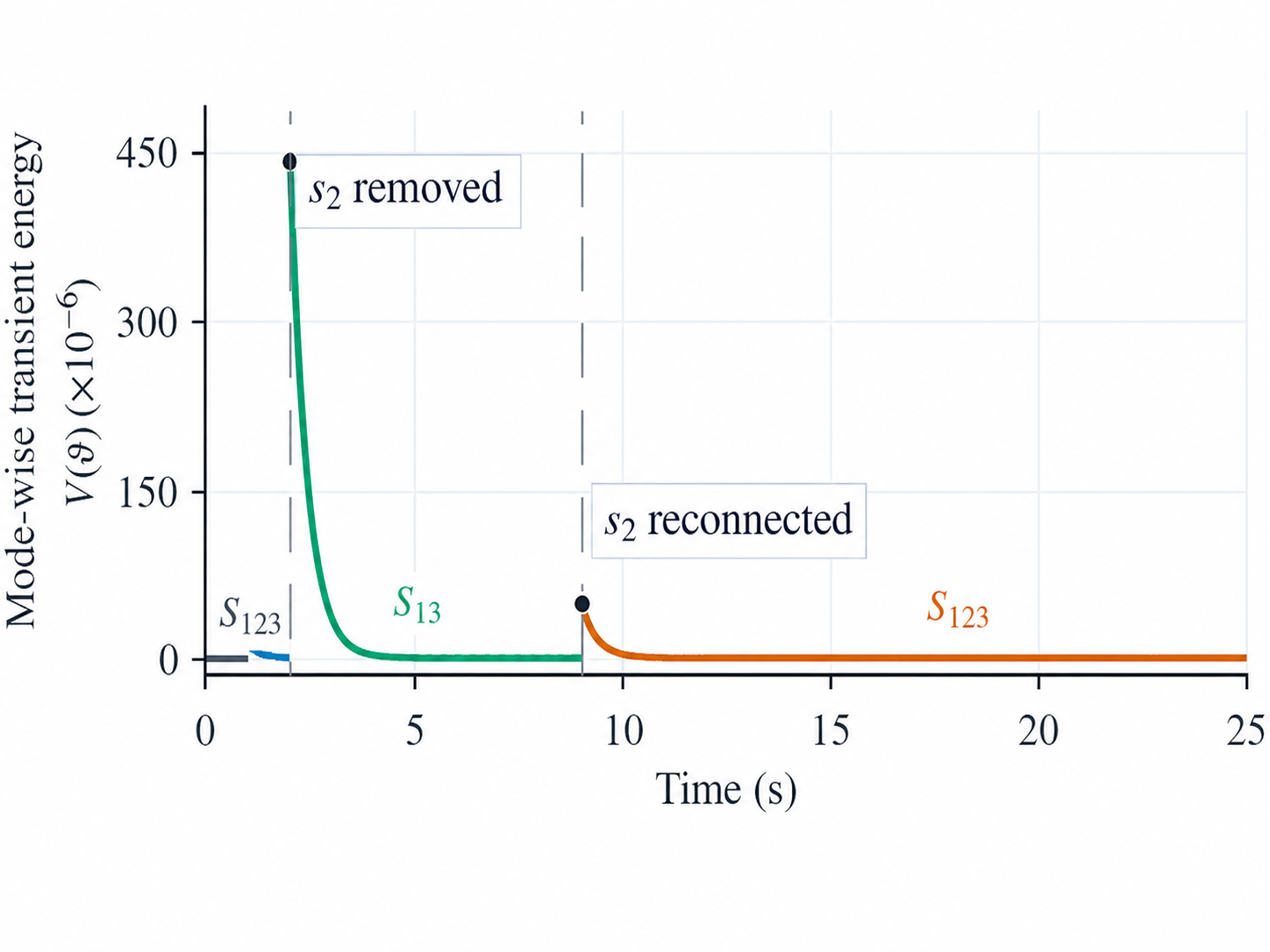}\\[-2pt]
 {\scriptsize (c) Mode-wise transient energy}
 \caption{Finite-schedule visualization for the illustrative example. The plots support the representative construction and dwell-time interpretation, but they do not replace the quotient-metric calculations.}
 \label{fig:ill_visual_summary}
\end{figure}

Figure~\ref{fig:ill_visual_summary} gives a compact visualization of the finite schedule. The frequency plot records the event epochs and the COI-relative speed responses; the reduced trajectory plot shows the switching path in low-dimensional internal coordinates; and the energy trace shows the mode-wise transient energy behavior, used as a local consistency check. These plots are explanatory; the quantitative cross-dimensional claims are the metric values in Table~\ref{tab:ill_event_compare}. The mode-wise energy transition and consistency details are provided in \ref{app:energy_consistency_details}.

The numerical energy check in Table~\ref{tab:ill_energy_consistency} is used only as a finite-schedule consistency check. It shows that, for the selected residence times, the realized post-switch states remain within the sampled local energy neighborhoods of the subsequent active modes. This should not be read as a computation of the minimum dwell time, nor as a full verification of a general switching-stability theorem for all admissible switching signals.

\begin{table}[!t]
 \centering
 \caption{Finite-schedule energy consistency check under the selected dwell time.}
 \label{tab:ill_energy_consistency}
 \scriptsize
 \begin{tabular}{@{}lcccc@{}}
 \toprule
 Event & Target & $V^+(\times 10^{-4})$ & $r_{S^+}(\times 10^{-4})$ & $r_{S^+}/V^+$\\
 \midrule
 Removal & $S_{13}$ & $4.4051$ & $6.6077$ & $1.50$ \\
 Reconnection & $S_{123}$ & $0.4544$ & $7.7637$ & $17.09$ \\
 \bottomrule
 \end{tabular}
\end{table}

Here $r_{S^+}$ is a sampled boundary value used as a numerical surrogate for the local energy-neighborhood radius in the target mode. A more valuable but substantially harder problem is to compute an explicit minimum admissible dwell time from the quotient-space switching gap, the affine offset, the mode-wise decay rate, and the reconnection policy. This quantitative dwell-time synthesis problem is left for future work.

This example therefore illustrates the role of the dwell-time condition in the quotient-space framework: the metric $d_{\mathcal V}$ quantifies the dimension-changing switching gap, while a sufficiently long residence interval allows the active mode to dissipate the induced deviation before the next structural event.

\section{Concluding Remarks}

The main contribution of this paper is to provide a universal state space for dimension-varying (control) systems. As mentioned in the Introduction, a distance between Euclidean spaces of different dimensions is needed for comparing trajectories that evolve in spaces of different dimensions. Different choices of such a distance are possible. This paper proposes a unified distance and provides evidence for its usefulness.

Under the unified distance space $\Omega$, dimension-varying (control) systems can be modeled in a common state space, and several classical problems for fixed-dimensional control systems can be formulated in this framework. Controllability, observability, stabilizability, and the disturbance decoupling problem have been discussed. Other control problems may also be formulated within this state space. The illustrative power-system example shows how the abstract quotient-space objects can be instantiated through translated reduced representatives, projection/lift benchmarks, event-wise $d_{\mathcal V}$-gap calculations, and a dwell-time consistency check.

An important remaining question is whether the distance is physically meaningful for a given application domain. As mentioned above, $d_{{\cal V}}(x,y)=0$ (i.e., $x\lra y$) is equivalent to the existence of $\J_{\a}$ and $\J_{\b}$ such that
$$
x\otimes \J_{\a}=y\otimes \J_{\b}.
$$
From an information-theoretic viewpoint, this equivalence indicates that $x$ and $y$ share the same replicated information pattern. For signal-processing systems, such a distance is natural. For mechanical systems and other physical domains, further practical verification remains necessary.

\appendix

\section{Power-System Model Reduction Details}
\label{app:power_model_details}

For each $S\in\mathcal S$, let $x^{(S)}=[\theta_i,\omega_i]_{i\in S}^{\top}$ denote the dynamic state and let $z^{(S)}$ collect the algebraic network variables. A mode-wise DAE representative is
\begin{align*}
 \dot\theta_i &= \omega_s\omega_i, \qquad i\in S,\\
 2H_i\dot\omega_i &= P_{m,i}-D_i\omega_i-P_{e,i}^{(S)}(x^{(S)},z^{(S)}), \qquad i\in S,\\
 0 &= g^{(S)}(x^{(S)},z^{(S)};p^{(S)}).
\end{align*}
When the algebraic Jacobian $\partial g^{(S)}/\partial z^{(S)}$ is nonsingular near the segment-local operating point, the implicit function theorem gives a local map $z^{(S)}=h^{(S)}(x^{(S)})$. Substitution yields a mode-wise ODE representative $\dot x^{(S)}=f^{(S)}(x^{(S)})$.

In the classical synchronous-generator model used in the example, the algebraic elimination is represented by the Kron-reduced internal network
\begin{equation*}
 I^{(S)}=Y_{\rm red}^{(S)}E^{(S)},
\end{equation*}
where $Y_{{\rm red},ij}^{(S)}=G_{{\rm red},ij}^{(S)}+jB_{{\rm red},ij}^{(S)}$ and $E_i=|E_i|e^{j\theta_i}$. The active- and reactive-power injections are
\begin{align*}
 P_{e,i}^{(S)}
 &=\sum_{j\in S}|E_i||E_j|
 \bigl(G_{{\rm red},ij}^{(S)}\cos\theta_{ij}
 +B_{{\rm red},ij}^{(S)}\sin\theta_{ij}\bigr),\\
 Q_{e,i}^{(S)}
 &=\sum_{j\in S}|E_i||E_j|
 \bigl(G_{{\rm red},ij}^{(S)}\sin\theta_{ij}
 -B_{{\rm red},ij}^{(S)}\cos\theta_{ij}\bigr),
\end{align*}
where $\theta_{ij}=\theta_i-\theta_j$. The swing equation uses the active-power balance, while the reactive-power relation determines the AC operating point used to construct the segment-local equilibrium. The transient-energy argument used in the example is interpreted within the standard local-potential representation of the classical direct method; for a lossy reduced network, nonconservative terms must be neglected or treated separately.

\section{Reduced Coordinates and Segment-Local Equilibria}
\label{app:reduced_coordinates_details}

For each $S\in\mathcal S$, the center-of-inertia frequency is
\begin{equation*}
 \omega_{\rm COI}^{(S)}=\frac{\sum_{i\in S}H_i\omega_i}{\sum_{i\in S}H_i}.
\end{equation*}
The example uses $s_1$ as the local reference machine and defines
\begin{equation*}
 \vartheta_{i1}^{(S)}={\rm wrap}_{(-\pi,\pi]}(\theta_i-\theta_1),
 \qquad
 \varpi_{i1}^{(S)}=\omega_i-\omega_1,
 \quad i\in S\setminus\{s_1\}.
\end{equation*}
The reference-machine basis and the COI-free basis span the same local internal subspace; they differ only by an invertible local coordinate change near the operating point. The wrapped angle coordinate is used only inside a local chart that does not cross the $\pm\pi$ branch cut.

The unshifted reduced representatives are
\begin{equation*}
 c^{(S_{123})}=[\vartheta_{21},\vartheta_{31},\varpi_{21},\varpi_{31}]^{\top}\in\mathbb R^4,
 \qquad
 c^{(S_{13})}=[\vartheta_{31},\varpi_{31}]^{\top}\in\mathbb R^2.
\end{equation*}
Table~\ref{tab:app_segment_equilibria} gives the segment-local equilibria used in the translated quotient-space analysis.

\begin{table}[!t]
    \centering
    \caption{Segment-local equilibria used for the translated quotient-space analysis.}
    \label{tab:app_segment_equilibria}
    \tiny
    \begingroup
    \setlength{\tabcolsep}{1.5pt}
    \renewcommand{\arraystretch}{1.08}
    \begin{tabular}{>{\raggedright\arraybackslash}p{0.16\columnwidth}
                    >{\raggedright\arraybackslash}p{0.38\columnwidth}
                    >{\raggedright\arraybackslash}p{0.29\columnwidth}
                    >{\centering\arraybackslash}p{0.08\columnwidth}}
        \toprule
        Mode & Raw synchronous operating point & Reduced equilibrium $c_*^{(S)}$ & Object \\
        \midrule
        $[1,2]$, $S_{123}$ &
        \begin{tabular}[t]{@{}l@{}}$(0.0795,-0.0010,-0.0642,$\\$-0.0010,0.0415,-0.0010)$\end{tabular} &
        $[-0.1437,-0.0380,0,0]^\top$ & $\bar 0$ \\
        $[2,9]$, $S_{13}$ &
        \begin{tabular}[t]{@{}l@{}}$(0.0795,-0.0078,$\\$0.0360,-0.0078)$\end{tabular} &
        $[-0.0435,0]^\top$ & $\bar 0$ \\
        $[9,25]$, $S_{123}$ &
        \begin{tabular}[t]{@{}l@{}}$(0.0795,-0.0058,-0.1016,$\\$-0.0058,0.0376,-0.0058)$\end{tabular} &
        $[-0.1811,-0.0419,0,0]^\top$ & $\bar 0$ \\
        \bottomrule
    \end{tabular}
    \endgroup
\end{table}

\section{Event-State Provenance and Metric Recalculation}
\label{app:event_metric_recalculation}

At $t=2.0\,{\rm s}$, removal of $s_2$ maps the four-dimensional three-machine representative to the two-dimensional two-machine representative. Immediately before removal,
\begin{equation*}
 c_{(123)}^-=
 \begin{bmatrix}-0.143615\\-0.037693\\4.03\times10^{-5}\\-1.28\times10^{-6}\end{bmatrix},
 \qquad
 \xi_{(123)}^-=
 \begin{bmatrix}6.36\times10^{-5}\\2.81\times10^{-4}\\4.03\times10^{-5}\\-1.28\times10^{-6}\end{bmatrix}.
\end{equation*}
Immediately after removal,
\begin{equation*}
 c_{(13)}^+=
 \begin{bmatrix}-0.037693\\-1.28\times10^{-6}\end{bmatrix},
 \qquad
 \xi_{(13)}^+=
 \begin{bmatrix}5.8361\times10^{-3}\\-1.28\times10^{-6}\end{bmatrix}.
\end{equation*}
The unshifted deletion law is
\begin{equation*}
 W_{\rm rem}^{\rm red}=\begin{bmatrix}0&1&0&0\\0&0&0&1\end{bmatrix},
 \qquad
 c_{(13)}^+=W_{\rm rem}^{\rm red}c_{(123)}^-.
\end{equation*}
In translated coordinates,
\begin{equation*}
 \xi_{(13)}^+=W_{\rm rem}^{\rm red}\xi_{(123)}^-+\eta_{\rm rem},
 \qquad
 \eta_{\rm rem}=\begin{bmatrix}5.55\times10^{-3}\\0\end{bmatrix}.
\end{equation*}
Using \eqref{eq:ill_dv_rule},
\begin{align*}
 \delta_{\rm rem}^{\rm shift}&=d_{\mathcal V}(\xi_{(123)}^-,\xi_{(13)}^+)=4.0056\times10^{-3},\\
 \delta_{\rm rem}^{\rm proj}&=d_{\mathcal V}(\xi_{(123)}^-,\Pi_2^4\xi_{(123)}^-)=7.83\times10^{-5},\\
 \delta_{\rm rem}^{\rm policy}&=d_{\mathcal V}(\xi_{(13)}^+,\Pi_2^4\xi_{(123)}^-)=4.0049\times10^{-3}.
\end{align*}

At $t=9.0\,{\rm s}$, $s_2$ is reinserted by the scenario-level policy labeled \texttt{synchronized\_zero\_power}. Immediately before reconnection,
\begin{equation*}
 c_{(13)}^-=
 \begin{bmatrix}-0.043534\\1.25\times10^{-6}\end{bmatrix},
 \qquad
 \xi_{(13)}^-=
 \begin{bmatrix}-5.73\times10^{-6}\\1.25\times10^{-6}\end{bmatrix}.
\end{equation*}
Immediately after reconnection,
\begin{equation*}
 c_{(123)}^+=
 \begin{bmatrix}-0.196802\\-0.043534\\5.55\times10^{-7}\\1.25\times10^{-6}\end{bmatrix},
 \qquad
 \xi_{(123)}^+=
 \begin{bmatrix}-1.5699\times10^{-2}\\-1.6024\times10^{-3}\\5.55\times10^{-7}\\1.25\times10^{-6}\end{bmatrix}.
\end{equation*}
The realized reduced reinsertion law is affine:
\begin{equation*}
 c^+=A_{\rm rec}c^-+b_{\rm rec},
 \quad
 A_{\rm rec}=\begin{bmatrix}0&0\\1&0\\0&0\\0&1\end{bmatrix},
 \quad
 b_{\rm rec}=\begin{bmatrix}\eta_{\vartheta}\\0\\\eta_{\varpi}\\0\end{bmatrix},
\end{equation*}
where $\eta_{\vartheta}=-0.196802$ and $\eta_{\varpi}=5.55\times10^{-7}$ for the analyzed event. The exact lift has zero intrinsic projection cost:
\begin{equation*}
 d_{\mathcal V}(\xi_{(13)}^-,\Pi_4^2\xi_{(13)}^-)=0.
\end{equation*}
The realized shifted and policy-induced gaps are
\begin{align*}
 \delta_{\rm rec}^{\rm shift}&=d_{\mathcal V}(\xi_{(13)}^-,\xi_{(123)}^+)=7.8869\times10^{-3},\\
 \delta_{\rm rec}^{\rm proj}&=d_{\mathcal V}(\xi_{(13)}^-,\Pi_4^2\xi_{(13)}^-)=0,\\
 \delta_{\rm rec}^{\rm policy}&=d_{\mathcal V}(\xi_{(123)}^+,\Pi_4^2\xi_{(13)}^-)=7.8869\times10^{-3}.
\end{align*}

The linear parts of the realized event maps have finite induced gains in the quotient metric:
\begin{align*}
 \|W_{\rm rem}^{\rm red}\|_{\mathcal V}
 &=\sqrt{\frac{4}{2}\lambda_{\max}\!\left((W_{\rm rem}^{\rm red})^{\top}W_{\rm rem}^{\rm red}\right)}=\sqrt2,\\
 \|A_{\rm rec}\|_{\mathcal V}
 &=\sqrt{\frac{2}{4}\lambda_{\max}\!\left(A_{\rm rec}^{\top}A_{\rm rec}\right)}=\frac{1}{\sqrt2}.
\end{align*}
The affine offsets $\eta_{\rm rem}$ and $\eta_{\rm rec}$ do not affect these Lipschitz constants, but they determine the post-event displacement from the common translated equilibrium object.

The two-machine translated vector field can also be lifted to the four-dimensional comparison space by
\begin{equation*}
 G_{13}(y)=g_{13}(\Pi_2^4y)\otimes\mathbf 1_2,
 \qquad y\in\Omega^4.
\end{equation*}
This projection-induced fixed-dimensional view illustrates the equivalent-switching-system construction. Note that it is only a quotient-compatible comparison view, not a claim that the disconnected physical mode is a four-dimensional autonomous power-system model.

\section{Mode-Wise Energy Consistency Details}
\label{app:energy_consistency_details}

\begin{figure}[!t]
 \centering
 \includegraphics[width=0.8\columnwidth]{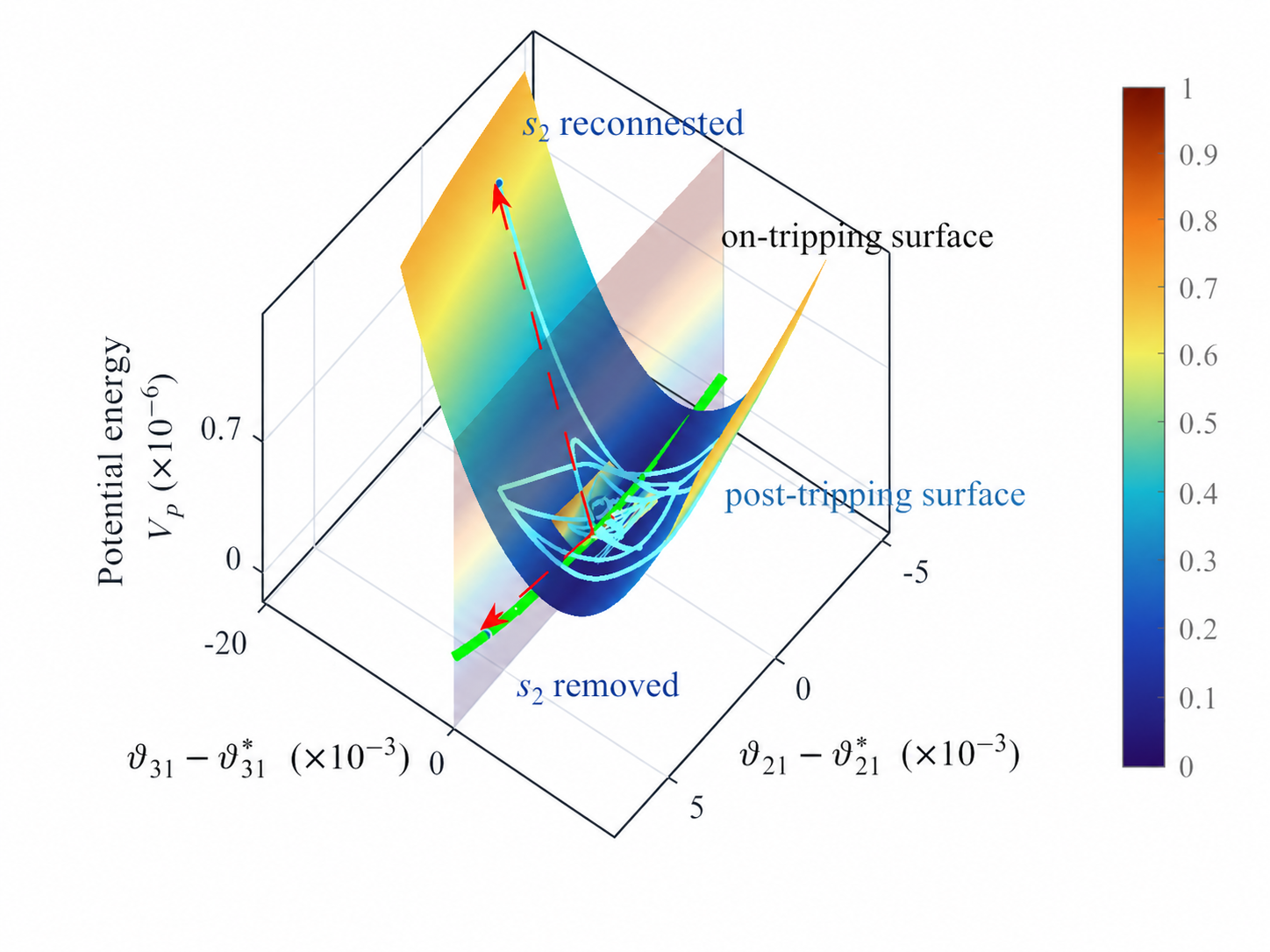}
 \caption{Energy-surface view of the translated reduced trajectory. This surface visualizes the mode-wise transient-energy landscape; it is not a $d_{\mathcal V}$ surface.}
 \label{fig:app_energy_surface}
\end{figure}

For a fixed $S\in\mathcal S$, the translated reduced swing dynamics can be written locally in the internal form
\begin{equation*}
 \dot\vartheta_\rho=\nu,
 \qquad
 \mathsf M_S\dot\nu=-\mathsf D_S\nu-\nabla_{\vartheta_\rho}\Delta U^{(S)}(\vartheta_\rho),
\end{equation*}
where $\mathsf M_S\succ0$ and the reduced damping is assumed to be effective on the internal frequency subspace. Define
\begin{align*}
 \Delta U^{(S)}(\vartheta_\rho)
 ={}&U^{(S)}(\tilde\vartheta_*^{(S)}+\vartheta_\rho)-U^{(S)}(\tilde\vartheta_*^{(S)})
 -\nabla U^{(S)}(\tilde\vartheta_*^{(S)})^{\top}\vartheta_\rho,\\
 V^{(S)}(\xi)
 ={}&\frac12\nu^{\top}\mathsf M_S\nu+\Delta U^{(S)}(\vartheta_\rho).
\end{align*}
Along a fixed active mode,
\begin{equation*}
 \dot V^{(S)}=-\nu^{\top}\mathsf D_S\nu\le 0.
\end{equation*}
This mode-wise energy argument is used only as a local consistency check for the selected finite schedule. It is not a common Lyapunov function on $\Omega^4\cup\Omega^2$ and does not by itself establish a general switching-stability theorem. Figure~\ref{fig:app_energy_surface} shows the energy surface visualization for the trajectory used in the illustrative example.

\balance

%
%

\end{document}